\documentstyle[12pt]{article}
\parskip=\medskipamount                 
\thicklines                     
\newcommand{\bim}[6]{\bibitem{#1}#2, {\em #3\/}$\;${\bf
#4}$\;$(#5)$\;${#6}.}
\def\inbar{\vrule height1.5ex width.4pt depth0pt}
\def\IN{\relax{\rm I\kern-.18em N}}
\def\IQ{\relax\,\hbox{$\inbar\kern-.3em{\rm Q}$}}
\def\IR{\relax{\rm I\kern-.18em R}}
\def\ZZ{\relax{\sf Z\kern-.4em Z}}
\def\a{\alpha} \def\b{\beta}

   \def\cD{{\cal D}}
   
  \def\cL{{\cal L}}

\newtheorem{theorem}{Theorem}[section]
\newtheorem{proposition}{Proposition}[section]

\newtheorem{lemma}{Lemma}[section]

\catcode`\@=11

\newif\if@fewtab\@fewtabtrue

\catcode`\@=11

\newif\if@fewtab\@fewtabtrue

{\count255=\time\divide\count255 by 60
\xdef\hourmin{\number\count255}
\multiply\count255 by-60\advance\count255 by\time
\xdef\hourmin{\hourmin:\ifnum\count255<10 0\fi\the\count255}}
\def\ps@draft{\let\@mkboth\@gobbletwo
    \def\@oddhead{}
    \def\@oddfoot
       {\hbox to 7 cm{$\scriptstyle Draft\ version:\ \draftdate$
       \hfil}\hskip -7cm\hfil\rm\thepage \hfil}
    \def\@evenhead{}\let\@evenfoot\@oddfoot}

\def\ceqno{\global\@fewtabfalse
    \ifcase\@eqcnt \def\@tempa{& & &}\or \def\@tempa{& &}
      \or \def\@tempa{&}
      \or\def\@tempa{}\fi\@tempa
{\rm(\theequation)}}

\def\aeqno#1{\global\@fewtabfalse
    \ifcase\@eqcnt \def\@tempa{& & &}\or \def\@tempa{& &}
      \or \def\@tempa{&}
      \or\def\@tempa{}\fi\@tempa
{\rm(\theequation,#1)}}

\def\label#1{\ifnum\draftcontrol=1
 \global\def\draftnote{$\scriptstyle #1$}\fi
 \@bsphack\if@filesw {\let\thepage\relax
   \def\protect{\noexpand\noexpand\noexpand}%
\xdef\@gtempa{\write\@auxout{\string
      \newlabel{#1}{{\@currentlabel}{\thepage}}}}}\@gtempa
   \if@nobreak \ifvmode\nobreak\fi\fi\fi
  \@esphack}

\def\alabel#1#2{\label{#1}\global\@fewtabfalse
    \ifcase\@eqcnt \def\@tempa{& & &}\or \def\@tempa{& &}
      \or \def\@tempa{&}
      \or\def\@tempa{}\fi\@tempa
{\hbox to 3cm{\phantom{\rm(\theequation,#2)}
\draftnote \hfil}\hskip -3cm {\rm(\theequation,#2)}}}

\def\clabel#1{\label{#1}\global\@fewtabfalse
    \ifcase\@eqcnt \def\@tempa{& & &}\or \def\@tempa{& &}
      \or \def\@tempa{&}
      \or\def\@tempa{}\fi\@tempa
{\hbox to 3cm{\phantom{\rm(\theequation)}
\draftnote \hfil}\hskip -3cm{\rm(\theequation)}}}

\def\eqnarray{\def\draftnote{{}}\global\@fewtabtrue
\stepcounter{equation}\let\@currentlabel=\theequation
\global\@eqnswtrue
\global\@eqcnt\z@\tabskip\@centering\let\\=\@eqncr
$$\halign to \displaywidth\bgroup\@eqnsel\hskip\@centering\@eqcnt\z@
  $\displaystyle\tabskip\z@{##}$&\global\@eqcnt\@ne
  \hskip 1\arraycolsep \hfil${##}$\hfil
  &\global\@eqcnt\tw@ \hskip 1\arraycolsep
$\displaystyle\tabskip\z@{##}$
\hfil  \tabskip\@centering&\global\@eqcnt\thr@@\llap{##}\tabskip\z@
\cr}

\def\endeqnarray{\@@eqncr\egroup
      \global\advance\c@equation\m@ne$$\global\@ignoretrue}

\def\@eqnnum{\hbox to 3cm{\phantom{\rm(\theequation)} \draftnote
                         \hfil}\hskip -3cm {\rm(\theequation)}}

\def\@@eqncr{\let\@tempa\relax
    \ifcase\@eqcnt \def\@tempa{& & &}\or \def\@tempa{& &}
      \or \def\@tempa{&}
      \or\def\@tempa{}
\fi\@tempa
\if@eqnsw
\if@fewtab\@eqnnum\fi
\stepcounter{equation}\fi\global
\@eqnswtrue\global\@eqcnt\z@\global\@fewtabtrue\cr}

\def\draftcite#1{\ifnum\draftcontrol=1#1\else{}\fi}

\def\@lbibitem[#1]#2{\item{}\hskip -3cm \hbox to 2cm
{\hfil$\scriptstyle\draftcite{#2}$}\hskip
1cm[\@biblabel{#1}]\if@filesw
     {\def\protect##1{\string ##1\space}\immediate
      \write\@auxout{\string\bibcite{#2}{#1}}}\fi\ignorespaces}

\def\@bibitem#1{\item\hskip -3cm \hbox to 2cm
{\hfil $\scriptstyle\draftcite{#1}$}\hskip 1cm
\if@filesw \immediate\write\@auxout
       {\string\bibcite{#1}{\the\value{\@listctr}}}\fi\ignorespaces}

\def\nsection#1{\section{#1}\setcounter{equation}{0}}

\def\draftdate{\number\month/\number\day/\number\year\ \ \ \hourmin }

\global\def\draftcontrol{0}
\catcode`\@=12

\def\theequation{{\thesection.\arabic{equation}}}

\def\qq{\begin{eqnarray}}
\def\qqq{\end{eqnarray}}

\newlength{\shiftwidth}
\def\shift#1{&&\hbox to \shiftwidth{\hfill $\displaystyle#1$}}
\newlength{\sshiftwidth}
\def\sshift#1{\lefteqn{\hbox to
\sshiftwidth{\hfill$\displaystyle#1$}}}

\def\ie{{\it i.e.\ }}
\def\eg{{\it e.g.\ }}

\def\rhs{{\it r.h.s.\ }}
\def\lhs{{\it l.h.s.\ }}

\def\Rhs{RHS\ }
\def\ihs{$\ZZ$HS\ }
\def\asl{ASL\ }

\def\Tr{\mathop{{\rm Tr}}\nolimits}
\def\Vol{\mathop{{\rm Vol}}\nolimits}

\def\p{^{\prime}}

\def\prosign{\mathop{{\rm sign}}\nolimits}
\def\sign#1{\prosign\left(#1\right)}

\def\spint#1{\int\limits^{+\infty}_{\scriptstyle -\infty \atop [{#1} =
0]}}

\def\max{\mathop{{\rm max}}\nolimits}
\def\promod{\mathop{{\rm mod}}\nolimits}
\def\mod#1{\;(\promod #1)}
\def\mmax#1{\max\{#1\}}
\def\intg#1{ \left[ {#1} \right] }
\def\Leg#1{ \left( {#1 \over K} \right) }

\def\Scs{S_{\rm CS}}
\def\csc{S^{(c)}_{\rm CS}}
\def\snm{S_n(M)}
\def\snst{S_n(S^3)}
\def\scnm{S_n^{(c)}(M)}

\def\ordH{| H_1(M,\ZZ)|}

\def\aa{{\a_1,\ldots,\a_N}}
\def\daa{d\a_1 \cdots d\a_N}
\def\aam{\a_1^{2m_1} \cdots \a_N^{2m_N}}
\def\aap{\a_1 + {1\over \xp_1}, \ldots, \a_N + {1\over \xp_N}}
\def\aaps{\a_1 + \xp^*_1, \ldots, \a_N + \xp^*_N}
\def\paap{\left( \aap \right)}

\def\zmk{Z(M;k)}
\def\zsk{Z(S^3;k)}
\def\zmo{Z(M;1)}
\def\zso{Z(S^3;1)}
\def\zpmk{Z\p(M;k)}
\def\zcmk{Z^{(c)}(M;k)}
\def\ztrmk{Z^{({\rm tr})}(M;k)}
\def\ztr{ Z^{({\rm tr})} }

\def\Jalk{J_\aa(\cL;k)}
\def\Japlk{J_{\aap}(\cL;k)}
\def\Japslk{J_{\aaps}(\cL;k)}
\def\Jalkr{J^{({\rm res})}_\aa(\cL;k)}

\def\lcw{\lambda_{\rm CW}}
\def\lnm{\lambda_n(M)}

\def\sjN{\sum_{j=1}^N}
\def\pjN{\prod_{j=1}^N}
\def\snzi{\sum_{n=0}^\infty}
\def\snoi{\sum_{n=1}^\infty}
\def\smnf{\sum_{m\leq {3\over 4} n} }
\def\smmN{\sum_{\mm\geq 0 \atop m_1+\cdots+m_N = m}}
\def\smj{\sum_{\mu_j = \pm 1} }
\def\ssignp{ \sjN \sign{\xp_j} }

\def\sta{ \sum_{ {3 - K \over 2} < \ta \leq \kt } }

\def\smb{\sum_{1 \leq m_j \leq 2\mb_j + 1 \atop (\ojN)} }

\def\ojN{1\leq j\leq N}

\def\ppaj{\left( \pjN \a_j \right)}
\def\ppajp{\left( \pjN \left(\a_j+{1\over \xp_j} \right) \right)}
\def\ppajps{\left( \pjN \left(\a_j + \xp_j^* \right) \right)}

\def\sajk{ \sum_{0\leq \a_j \leq K-1 \atop (1\leq j \leq N)} }
\def\psum{\mathop{\sum\nolimits^{\prime}}\limits}
\def\spajk{ \psum \limits_{
  -K\leq \a_j \leq K \atop{ \a_j \in 2\ZZ+1
  \atop (1\leq j\leq N) } } }

\def\spak{ \psum \limits_{
  -K\leq \a \leq K \atop \a \in 2\ZZ+1 } }
\def\spgk{ \psum \limits_{
  -K\leq \gamma \leq K \atop \gamma \in 2\ZZ+1 } }
\def\spbk{ \psum \limits_{
  -K\leq \b \leq K \atop \b \in 2\ZZ+1 } }
\def\sak{ \sum_{\a=0}^{K-1} }

\def\iptk{{i\pi\over 2K}}
\def\pk{{\pi\over K}}
\def\ppk{\left(\pk \right)}
\def\knt{K^{- {N\over 2}}}
\def\tpik{{2\pi i\over K}}
\def\ptpik{\left(\tpik \right)}
\def\kt{{K-1\over 2}}
\def\kpt{{K+1\over 2}}
\def\ipk{{i\pi\over K}}
\def\pipk{\left( \ipk \right)}

\def\mm{m_1,\ldots,m_N}
\def\mmp{m_1\p,\ldots,m_N\p}
\def\Dmn{D_{m,n}}
\def\Dpmn{D^{(m,n)}_{\mm}}
\def\Dtmn{\tilde{D}_{m,n}}
\def\Dtpmn{\tilde{D}^{(m,n)}_{\mm}}
\def\cnmN{C^{(n)}_{\mm}}
\def\cnmNp{C^{(n)}_{\mmp}}
\def\Dtc{ \ppaj \smnf \Dtmn(\aa) }
\def\Dtcps{ \ppajps \smnf \Dtmn(\aaps) }
\def\Dtcp{ \ppajp \smnf \Dtmn(\aap) }

\def\polN{ P_{m_1}\left({\a_1 - 1\over 2}\right)
  \cdots P_{m_N}\left({\a_N - 1\over 2}\right) }
\def\polNp{ P_{m_1\p}\left({\a_1 - 1\over 2}\right)
  \cdots P_{m_N\p}\left({\a_N - 1\over 2}\right) }

\def\spint#1{\int\limits^{+\infty}_{\scriptstyle -\infty \atop [{#1} =
   0]}}

\def\gpq#1{ G(\xp,\xq;#1)}
\def\gpqi#1#2{ G_{#2}(\xp,\xq;#1)}

\def\pmap{P_m\left( {\a - 1 + {1\over \xp} \over 2} \right) }

\def\pmaps{P_m\left( {\a - 1 + \xp^* \over 2} \right) }

\def\pKap{P_K\left( {\a - 1 + {1\over \xp} \over 2} \right) }
\def\pKapm{P_K\left( {-\a - 1 + {1\over \xp} \over 2} \right) }
\def\pKaps{P_K\left( {\a - 1 + \xp^* \over 2} \right) }
\def\pKapsm{P_K\left( {-\a - 1 + \xp^* \over 2} \right) }
\def\pkev{P_K^{({\rm ev})}}

\def\eikap{\frac {e^ {i {\pi\over 4} (\kappa -1)}} {\sqrt{K}}}
\def\eikapw{ e^ {i {\pi\over 4} (\kappa -1) } }

\def\opr{${\mbox{Ohtsuki}}^\prime\;\;$}

\def\xp{p}
\def\xq{q}
\def\xr{r}
\def\xs{s}
\def\xx{x}
\def\spq{\check{q}}

\def\tU{\tilde{U}}
\def\tL{\tilde{L}}
\def\tcL{\tilde{\cL}}
\def\chU{\check{U}}
\def\chS{\check{S}}
\def\chT{\check{T}}

\def\mb{\bar{m}}
\def\tb{\bar{t}}
\def\ta{\tilde{\a}}
\def\zq{\ZZ[\spq]}
\def\zx{\ZZ[\xx]}
\def\zkx{\ZZ_K[\xx]}
\def\zpkx{\ZZ\p_K[\xx]}
\def\zrat{\ZZ\left[{1\over 2}, {1\over 3}, \ldots, {1\over 2n+1},
  {1\over \ordH} \right]}
\def\zratm{\ZZ\left[{1\over 2}, {1\over 3}, \ldots, {1\over 2n},
  {1\over \ordH} \right]}

\def\nz{N_{\rm zero}}

\def\seif{X\left( {p_1\over q_1}, \ldots, {p_n\over q_n} \right) }
\def\spk{ \sin \left( {\pi\over K} \right) }
\def\spbaj{ \sin \left( {\pi\over K} \b\a_j \right) }
\def\spb{ \sin \left( {\pi\over K} \b \right) }

\def\sph{ \sign{H\over P} }
\def\lcon{ L^{({\rm con }) } }
\def\spqj{ \sign{q_j\over p_j} }
\def\lenpqj{ L_{-\xp_j,q_j} }
\def\eqpb{ \spq^{-2^* p_j^* \b} - \spq^{2^* p_j^* \b} }
\def\eqh{ \spq^{-2^*} - \spq^{2^*} }
\def\eqhm{ \spq^{2^*} - \spq^{-2^*} }
\def\eqb{ \spq^{-2^*\b} - \spq^{2^*\b} }
\def\pps{ p_1^*,\ldots,p_N^* }
\def\qPHb{ \spq^{-4^* P^* H \b^2} }
\def\lqexp{ \spq^{ 4^* P^* H - 3\cdot 4^* \sph - 3\sjN
    s^\vee(q_j,p_j)} }
\def\intinfb{ \int_{-\infty}^{+\infty} d\b\, }
\def\lint{ \spint{\b} d\b\,
  e^{-{i\pi\over 2K} {H\over P} \b^2}
  \frac {\pjN \sin \left( \pk {\beta\over p_j} \right) }
  {\sin^{N-2} \left( \pk \b \right) } }

\begin{document}

\begin{titlepage}
\centerline{\hfill                 UMTG-183-95}
\centerline{\hfill                 q-alg/9504015}
\vfill
\begin{center}

{\large \bf
Witten's Invariants of Rational Homology Spheres at Prime Values of
$K$ and Trivial Connection Contribution.
} \\

\bigskip
\centerline{L. Rozansky\footnote{Work supported
by the National Science Foundation
under Grant No. PHY-92 09978.
}}

\centerline{\em Physics Department, University of Miami
}
\centerline{\em P. O. Box 248046, Coral Gables, FL 33124, U.S.A.}
\centerline{{\em E-mail address: rozansky@phyvax.ir.miami.edu}}

\vfill
{\bf Abstract}

\end{center}
\begin{quotation}

We establish a relation between the coefficients of asymptotic
expansion of the trivial connection contribution to Witten's
invariant of rational homology spheres and the invariants that
T.~Ohtsuki extracted from Witten's invariant at prime values of $K$.
We also rederive the properties of prime $K$ invariants discovered by
H.~Murakami and T.~Ohtsuki. We do this by using the bounds on Taylor
series expansion of the Jones polynomial of algebraically split
links, studied in our previous paper. These bounds are enough to
prove that Ohtsuki's invariants are of finite type.
The relation between Ohtsuki's invariants and trivial connection
contribution is verified explicitly for lens spaces and Seifert
manifolds.

\end{quotation}
\vfill
\end{titlepage}

\pagebreak
\nsection{Introduction}
\label{s1}
Witten's invariant of 3d manifolds defined in~\cite{Wi} by a path
integral over the $SU(2)$ connections $A_\mu$ on a 3d manifold $M$
\qq
&{\displaystyle
\zmk = \int [\cD A_\mu] e^{ {ik \over 2\pi} \Scs [A_\mu] },
}&
\label{1.1}
\\
&{\displaystyle
\Scs = {1\over 2} \Tr \epsilon^{\mu\nu\rho}
\int_M d^3 x\, \left( A_\mu \partial_\nu A_\rho +
{2\over 3} A_\mu A_\nu A_\rho \right)
}
&
\label{1.02}
\qqq
($k\in \ZZ$, $\Tr$ is the trace taken in the fundamental
representation of $SU(2)$) can also be calculated combinatorially
with the help of the surgery formula. Let $M$ be a 3d manifold
constructed by $(\xp_j,1)$ surgeries on the components $\cL_j$ of an
$N$-component link $\cL$ in $S^3$. A $(\xp,1)$ surgery means that the
meridian of the tubular neighborhood is glued to the parallel plus
$\xp$ meridians on the boundary of the knot complement (in other
words, a Dehn's surgery is performed on a knot with framing number
$\xp$). The invariant of $M$ can be expressed in terms of the framing
independent colored Jones polynomial $\Jalk$ of the link $\cL$:
\qq
\lefteqn{
{\zmk \over \zsk} =
\left( {2\over K} \right)^{N\over 2}
\exp \left[ {3\over 4} \pi i \left( {2\over K} - 1 \right)
\sjN \sign{\xp_j} \right]
}\hspace*{1in}
\label{1.2}
\\
&\displaystyle
\times
\sajk \Jalk \exp \left( \iptk \sjN \xp_j (\a_j^2 - 1) \right)
\pjN \sin \left( \pk \a_j \right),
&
\nonumber\\
&\displaystyle
\zsk = \sqrt{2\over K} \sin \pk, \qquad K = k + 2.
\label{1.3}
\qqq
The Jones polynomial $\Jalk$ is normalized in such a way that it is
multiplicative  for unlinked links and
$J(\mbox{empty link};k) = 1$,
$J_\a(\mbox{unknot}; k) = \frac
{\sin \left( \pk \a \right)} {\sin \left( \pk \right) }$.

Although N.~Reshetikhin and V.~Turaev proved~\cite{ReTu} that
eq.~(\ref{1.2}) indeed defines an invariant of $M$ (\ie the \lhs of
eq.~(\ref{1.2}) is invariant under Kirby moves), the topological
origin of this invariant remains somewhat obscure. The question is:
which of the ``classical'' topological invariants of $M$ are
contained inside $\zmk$? Two distinct approaches to this problem have
been tried. The first one is to study $\zmk$ for some particular
values of $K$. R.~Kirby and P.~Melvin discovered~\cite{KiMe} that if
$K$ is odd, then $\zmk$ is proportional to $\zmo$:
\qq
{\zmk\over \zsk} = \zpmk \times
\left\{
\begin{array}{cl}
{\zmo\over \zso}, &\mbox{if $K=-1 \mod{4}$}\\
{\overline{\zmo}\over\zso}, &\mbox{if $K=1 \mod{4}$},
\end{array}
\right.
\label{1.4}
\qqq
If $M$ is a rational homology sphere (\ihs), then $\zmo = \zso$ so
that $\zpmk = {\zmk \over \zsk}$.

The new invariant $\zpmk$ can be calculated by the following surgery
formula
\qq
\lefteqn{
\zpmk = \knt \exp \left[ - {\pi i\over 4}
\left( 3 + \kappa - {6\over K} \right) \sjN \sign{\xp_j} \right]
}\hspace*{1in}
\label{1.5}
\\
&&
\times
\spajk \Jalk \exp \left( \iptk \sjN \xp_j (\a_j^2 - 1) \right)
\pjN \sin \left(\pk\a_j\right),
\nonumber
\qqq
here
\qq
\kappa =
\left\{
\begin{array}{cl}
1 &\mbox{if $K=-1 \mod{4}$}\\
-1&\mbox{if $K=1 \mod{4}$},
\end{array}
\right.
\label{1.6}
\qqq
while $\psum$ means that we add an extra factor of ${1\over 2}$  to
the terms corresponding to the boundary values of summation index
($\a_j = \pm K$ in this case). We changed slightly the original
formula of~\cite{KiMe}: instead of taking a sum over $1\leq\a_j \leq
\kt$ we sum over odd $\a_j$ between 1 and $K-1$. This allows us to
get rid of some phase factors. We also double the range of summation
to $1-K\leq \a_j \leq K-1$ by using the fact that $\Jalk$ is an odd
function of its indices (we use the $2K$ periodicity in $\a_j$ in
order to extend $\Jalk$ to negative values of $\a_j$). Note that the
whole summand of eq.~(\ref{1.5}) has a periodicity of $2K$.

S.~Garoufalidis~\cite{Ga1} used nice properties of the gaussian sum
$\sak \exp\left( \tpik \a^2 \right)$ for prime values of $K$ in order
to study Witten's invariant of lens spaces and Seifert manifolds.
H.~Murakami and T.~Ohtsuki~\cite{Mu1}--\cite{Oh2} carried out a
detailed study of the invariant $\zpmk$ of rational homology
spheres (\Rhs) for prime $K$.
\begin{theorem}[H.~Murakami,~\cite{Mu1},~\cite{Mu2}]
\label{t1.1}
For a \Rhs $M$ and a prime $K>2$
\qq
\zpmk \in \zq, \qquad
\spq = e^{\tpik},
\label{1.7}
\qqq
$\zq$ being a cyclotomic ring.
\end{theorem}

We need more notations in order to present the results of Ohtsuki's
papers~\cite{Oh1},~\cite{Oh2}. We introduce a new variable
\qq
\xx = q - 1.
\label{1.8}
\qqq
A polynomial from $\zq$ can be reexpressed as a polynomial in $\xx$
with integer coefficients. It is defined modulo the polynomial
\qq
{(1 + \xx)^K - 1 \over \xx} =
\sum_{n=0}^{K-1} \frac {K (K-1) \cdots (K-n)} {(n+1)!} \xx^n,
\label{1.9}
\qqq
which is identically equal to zero for $\xx = e^{\tpik} - 1$. All the
coefficients of this polynomial except the one at $\xx^{K-1}$, are
divisible by $K$. As a result, all the coefficients at $\xx^n$,
$n\leq K-2$ for a polynomial of $\xx$ coming from $\zq$ are well
defined modulo $K$. We will limit our attention to the powers of
$\xx$ up to $\xx^{\kt}$. They are all well defined as elements of
$\ZZ_K$ if $K\geq 3$. Thus there is a homomorphism of rings:
\qq
^\diamond:\;\;
\zq \rightarrow \zpkx \stackrel{\rm def.}{=} \zkx/ \xx^{\kpt} \zkx.
\label{1.10}
\qqq
There is another homomorphism from polynomials (of maybe infinite
degree) of $\xx$ with rational coefficients to $\zpkx$:
\qq
^\vee:\;\; \IQ[[x]] \rightarrow \zpkx.
\label{1.11}
\qqq
The action of the operation $^\vee$ on rational numbers was
introduced in relation to Witten's invariants at prime values of $K$
by S.~Garoufalidis~\cite{Ga1}:
\qq
^\vee:\;\;\IQ\rightarrow \ZZ_K,\qquad
\left({\xp\over \xq}\right)^{\vee}= \xp\xq^*,
\label{1.12}
\qqq
here $\xq^*$ is the inverse of $\xq$ modulo $K$: $\xq\xq^* = 1
\mod{K}$.
The homomorphism
$^\vee$ acts on polynomials (infinite series)
by
removing all powers of $\xx$ higher than $\kt$
and
converting the remaining coefficients to $\ZZ_K$.

Now we can present (a slightly stronger version of) Ohtsuki's
results:
\begin{theorem}[T.~Ohtsuki~\cite{Oh1},~\cite{Oh2}]
\label{t1.2}
For any \Rhs $M$ there exists a sequence of rational numbers
$\lambda_n(M) \in \zratm \subset \IQ$, $n\geq 0$ so that for any
prime number $K$ such that
$\ordH\neq 0 \mod{K}$
\qq
\left[ \ordH \left({\ordH\over K}\right) \zpmk \right]^\diamond
= \left[ \snzi \lambda_n(M) \xx^n \right]^{\vee},
\label{1.13}
\qqq
here $\left( {\cdot\over K} \right)$ is the Legendre symbol.
\end{theorem}
We have slightly modified the theorem of ~\cite{Oh2}: Ohtsuki
required that $K>\ordH$,  he estimated that $\lnm \in \zrat \subset
\IQ$, $n>0$ and
he used $\zpkx = \zkx / \xx^{\kt} \zx$
instead of~(\ref{1.10}) (in other words, he did not fix the
coefficient at $\xx^{\kt}$). Murakami showed that
\qq
\lambda_0 = 1,\qquad \lambda_1 = \lcw,
\label{1.14}
\qqq
here $\lcw$ is the Casson-Walker invariant of \Rhs.

The second approach to the search of the topological meaning of
Witten's invariant $\zmk$ is based on the path integral
representation~(\ref{1.1}). According to quantum field theory, this
integral can be calculated by stationary phase approximation when
$K\rightarrow \infty$. The invariant is presented as a sum of
contributions coming from connected components $c$ of the moduli
space of flat connections on the manifold $M$:
\qq
\zmk = \sum_c \zcmk.
\label{1.15}
\qqq
Each contribution has a general form
\qq
\zcmk = \left( {4\pi^2\over K} \right)^{\nz\over 2}
{1\over \Vol(H_c)}
\exp \left[ 2\pi ik \csc + \snzi \scnm \pipk^n \right],
\label{1.16}
\qqq
here $H_c$ is the isotropy group, $\nz = \dim H_c^0 - \dim H_c^1$,
$H_c^{0,1}$ being the cohomologies of 0,1-forms taking values in the
adjoint $su(2)$ bundle, $\csc$ is the Chern-Simons action and $\snzi
\scnm \pipk^n$ is an asymptotic series. The coefficients $\scnm$ are
called $(n+1)$-loop corrections. They might be related to
``classical'' topological invariants of $M$. Indeed, the 1-loop
correction $S_0^{(c)}$ is related to the Reidemeister-Ray-Singer
torsion. An attempt to relate the asymptotic properties of the
surgery formula~(\ref{1.2}) for lens spaces and Seifert manifolds to
the quantum field theory predictions of
eqs.~(\ref{1.15}),~(\ref{1.16}) was initiated by D.~Freed and
R.~Gompf~\cite{FrGo} and carried out further by L.~Jeffrey~\cite{Je},
S.~Garoufalidis~\cite{Ga1} and also in the papers~\cite{Ro*1,Ro*2}. A
complete agreement between the surgery formula and 1-loop predictions
was observed.

If the manifold $M$ is a \Rhs, then the trivial connection is a
separate point in the moduli space of flat connections. According to
quantum field theory, its contribution is of the form
\qq
\ztrmk = {\sqrt{2} \pi \over K^{3\over 2} \ordH^{3\over 2} }
\exp \left[ \snoi \snm \pipk^n \right].
\label{1.17}
\qqq
A representation of the coefficients $\snm$  in terms of $(n+1)$-loop
Feynman diagrams was carried out by S.~Axelrod and
I.~Singer~\cite{AxSi}, M.~Kontsevich~\cite{Ko}, C.~Taubes~\cite{Ta}
and others. We studied how
the trivial connection contribution can be extracted
from the surgery formula~(\ref{1.2}). We derived a knot surgery
formula~\cite{Ro1} and a link surgery formula~\cite{Ro3} for it. The
knot formula allowed us to show that
\qq
S_1(M) = 6\lcw.
\label{1.017}
\qqq

The link surgery formula of~\cite{Ro3} was much less explicit than
the knot formula of~\cite{Ro1}, because it did not express $\ztrmk$
directly in terms of derivatives of the Jones polynomial $\Jalk$.
However we derived an explicit surgery formula~\cite{Ro4} for
algebraically split links (\asl).

In this paper we are going to prove the following:
\begin{proposition}
\label{p1.1}
Ohtsuki's invariants $\lnm$ of eq.~(\ref{1.13}) and loop corrections
to the trivial connection contribution $\snm$ of eq.~(\ref{1.17}) can
be expressed in terms of each other through the following relation
\qq
\hspace*{-0.1in}
\snzi \lnm \xx^n = \exp \left[ \snoi (\snm - \snst) \pipk^n \right]
\equiv {\ppk \over \sin \ppk}
\exp \left[ \snoi \snm \pipk^n \right]
\label{1.18}
\qqq
by substituting either
\qq
\xx = e^{\tpik} - 1 =
\tpik \snzi {1\over (n+1)!} \ptpik^n,
\label{1.19}
\qqq
or
\qq
\ipk = {1\over 2} \log(1+\xx) =
{\xx\over 2} \snzi (-1)^n {\xx^n\over n+1}.
\label{1.20}
\qqq
\end{proposition}
In other words, we will show that
\qq
\left[ \ordH \left( {\ordH\over K} \right)
\zpmk \right]^{\diamond} =
\left[ \ordH^{3\over 2} {\ztrmk\over \zsk} \right]^\vee.
\label{1.21}
\qqq
In process of doing this we will rederive the results
of~\cite{Mu1}-\cite{Oh2} maybe in a somewhat more explicit way.

Our proof of the Proposition~\ref{p1.1}
will be based on the following
two propositions derived in ~\cite{Ro2} and~\cite{Ro4} by using some
``physical'' considerations:
\begin{proposition}
\label{p1.2}
Let $\cL$ be an algebraically split link (\asl) in $S^3$. Then its
framing-independent colored Jones polynomial has the
following Taylor series expansion in powers of $K$:
\qq
\Jalk = \ppaj \snzi \smnf \Dmn(\aa) \pipk^n,
\label{1.22}
\qqq
here $\Dmn(\aa)$ are even homogeneous polynomials of degree $2m$:
\qq
\Dmn(\aa) = \smmN \Dpmn \aam
\label{1.23}
\qqq
and
\qq
m_j \leq n-m, \qquad 1\leq j\leq N.
\label{1.24}
\qqq
\end{proposition}
\begin{proposition}
\label{p1.3}
Let $M$ be a rational homology sphere (\Rhs) constructed by
$(\xp_j,1)$ surgeries on the components of an $N$-component link
$\cL$ in $S^3$. Then the loop corrections to the trivial connection
contribution~(\ref{1.17}) to Witten's invariant of $M$ are given by
the formula:
\qq
\lefteqn{
\exp \left[ \snoi (\snm - \snst) \pipk^n \right]
}
\label{1.25}
\\
& &=
\exp \left[ {3\over 4} \pi i \left( {2\over K} - 1 \right)
\sjN \sign{\xp_j} \right]
\exp \left[ - \iptk \sjN \left( \xp_j + {1\over p_j} \right) \right]
{i^N\over (2K)^{N\over 2} }
\left| \pjN \xp_j \right|^{3\over 2}
\nonumber\\
&&
\qquad
\times
\spint{\a_j}
\daa\, \exp \left( \iptk \sjN \xp_j\a_j^2 \right)
\Japlk.
\nonumber
\qqq
The symbol $\spint{\a_j}$ means that the integral has to be
calculated in the following way: first, an expansion~(\ref{1.22}) has
to be substituted and then the gaussian integrals over $\a_j$ have to
be calculated for each polynomial $\ppajp \Dmn\paap$ separately (for
more details see~\cite{Ro4}).
\end{proposition}
The Proposition~\ref{p1.2} is essential for all our calculations. The
Proposition~\ref{1.3} is needed only for the derivation of
eq.~(\ref{1.21}). In other words, we could use eq.~(\ref{1.25}) as a
definition of $\snm$ and then prove eq.~(\ref{1.18}) which amounts to
proving the Theorem~\ref{t1.2}.

In Section~\ref{s2} we modify the surgery formula~(\ref{1.5}) and
prove the Theorem~\ref{t1.1}. Our main tool is the observation that
the gaussian sum $\sak \spq^{\a^2}$ is proportional to $\xx^{\kt}$,
while the sum $\sak \spq^{\a^2} \a^{2m}$ for $m\leq \kt$ is only
proportional to $\xx^{\kt-m}$. This is similar to the behavior of
gaussian integrals: each two extra powers of $\a$ in the integral
$\int_{-\infty}^{+\infty} e^{\tpik\a^2}\a^{2m}\, d\a$ bring a power of
$K$ to denominator. In
Section~\ref{s3}
we prove the
Proposition~\ref{p1.1} and thus also the Theorem~\ref{t1.2}. We use
again the similarities of the formulas for
$\sak \spq^{\a^2} \a^{2m}$ and
$\int_{-\infty}^{+\infty} e^{\tpik\a^2}\a^{2m}\, d\a$. These
similarities are due to the fact that both the sum
$\sak \spq^{\a^2 + 2n\a}$ and the integral
$\int_{-\infty}^{+\infty} e^{\tpik (\a^2 + 2n\a)}\,d\a$
can be calculated by completing the square in the exponents. In
Section~\ref{s4}
we derive a rational surgery formula for
$\zpmk$ which is similar to the
formula~(\ref{4.1}) of~\cite{Je}
for the
original Witten's invariant $\zmk$.
We use this formula to verify the Proposition~\ref{p1.1} for lens
spaces and Seifert manifolds which are rational homology spheres.
In
Section~\ref{s5}
we discuss
the properties of Ohtsuki's invariants $\lnm$ as related to the
properties of invariants $\snm$ studied in~\cite{Ro4}.

\nsection{Gaussian Sums and Divisibility in Cyclotomic Ring}
\label{s2}
We start by modifying the surgery formula~(\ref{1.5}). Since
${1\over 4}(\a_j^2 - 1) \in \ZZ$,
then
\qq
\exp \left[ \iptk \sjN \xp_j (\a_j^2 - 1) \right]
= \spq^{{1\over 4} \sjN \xp_j (\a_j^2 - 1)}
= \spq^{4^*\sjN \xp_j (\a_j^2 - 1)}.
\label{2.1}
\qqq
Also since
${1\over 2}(\a_j \pm \sign{\xp_j}) \in \ZZ$,
\qq
\sin \left( \pk \a_j \right) =
{i\over 2} \smj \mu_j e^{-\ipk \mu_j \a_j} =
{i\over 2} \smj \mu_j \spq^{-2^* \mu_j\a_j}
\spq^{\left( 2^* - {1\over 2} \right) \sign{\xp_j}}.
\label{2.2}
\qqq

The Jones polynomial $\Jalk$ is odd and
$e^{\iptk \xp_j \a_j^2}$
is even as a function of $\a_j$. Therefore we can drop the factor
${1\over 2}$ and put $\mu_j = 1$ in eq.~(\ref{2.2}) upon substituting
it into eq.~(\ref{1.5}):
\qq
\zpmk & = &
\knt i^N
e^{-{i\pi\over 4} (3 + \kappa) \ssignp}
\spq^{ \left( {1\over 4} + 2^* \right) \ssignp}
\spq^{-4^*\sjN \xp_j}
\label{2.3}
\\
&&
\qquad
\times
\spajk
\spq^{\sjN ( 4^* \xp_j \a_j^2 - 2^* \a_j ) }
\Jalk.
\nonumber
\qqq
After completing the square
\qq
4^* \xp_j \a_j^2 - 2^* \a_j =
4^* \xp_j (\a_j - \xp_j^*)^2 - 4^* \xp_j^* \mod{K},
\qqq
we shift the summation variable $\a_j$ by $\xp_j^*$ (we assume that
$\xp_j^*$ is even in order to preserve the parity of $\a_j$, we can
always make such choice of $\xp_j^*$ since $K$ is odd). Then
\qq
\zpmk & = &
\knt i^N
e^{-{i\pi\over 4} (3 + \kappa) \ssignp}
\spq^{ \left( {1\over 4} + 2^* \right) \ssignp}
\spq^{-4^* \sjN(\xp_j + \xp_j^*) }
\label{2.4}
\\
&&
\qquad
\times
\spajk
\spq^{ 4^* \sjN \xp_j \a_j^2}
\Japslk.
\nonumber
\qqq
Next we use the identities
\qq
&\displaystyle
ie^{-{i\pi\over 2} \sign{\xp_j}} = \sign{\xp_j},
\label{2.5}
\\
&\displaystyle
{1 - \kappa K \over 4} = 4^*
\label{2.6}
\qqq
in order to rearrange the phase factors preceding the sum in
eq.~(\ref{2.4}):
\qq
\zpmk & = & \knt e^{ {i\pi\over 4} (\kappa - 1) \ssignp }
\left( \pjN \sign{\xp_j} \right)
\spq^{ 4^* \sjN ( 3\sign{\xp_j} - \xp_j - \xp_j^*) }
\nonumber\\
&&\times
\spajk \spq^{ 4^* \sjN \xp_j \a_j^2 }
\Japslk.
\label{2.7}
\qqq

Since the values of the Jones polynomial $\Jalk$ belong to the
cyclotomic ring $\zq$ when all the indices $\a_j$ are odd, we can
apply to it a combination of Proposition~\ref{p1.2} and Lemma 2.3
of~\cite{Mu1}:
\begin{proposition}
\label{p2.1}
For odd values of its indices $\a_j$, the unframed colored Jones
polynomial of an $N$-component \asl $\cL$ in $S^3$ can be presented
as the following sum:
\qq
\Jalk & = &
\ppaj \sum_{n=0}^{ (N+1)\kt } \smnf \Dtmn (\aa) \xx^n
\label{2.8}
\\
\qquad &&
+ \xx^{ (N+1)\kt + 1 } \Jalkr,
\qquad \Jalkr \in \zq,
\nonumber
\qqq
here $\Dtmn(\aa)$ are homogeneous polynomials of degree $2m$:
\qq
&\displaystyle
\Dtmn (\aa) =
\smmN \Dtpmn \aam,
\label{2.9}
\\
&\displaystyle
m_j \leq n - m,
\label{2.10}
\qqq
and the polynomials $\ppaj \smnf \Dtmn (\aa)$ take integer values
when $\a_j$ are odd.
\end{proposition}
The latter property of the polynomials $\Dtmn$ allows us to express
them in terms of ``binomial coefficient'' polynomials:
\qq
&\displaystyle
\Dtc = \smb \cnmN \polN,
\nonumber\\
&\displaystyle
\cnmN \in \ZZ,
\label{2.11}
\qqq
here
\qq
P_m(\a) = \frac {\a (\a - 1) \cdots (\a - m + 1)} {m!}
\label{2.12}
\qqq
and $\mb_j$ are the maximum values of $m_j$ in the
representation~(\ref{2.9}) of all the polynomials $\Dtmn$ appearing
in the \lhs of eq.~(\ref{2.11}).

The following proposition is a reflection of the
inequality~(\ref{2.10}) for the representation~(\ref{2.11}):
\begin{proposition}
\label{pi2.1}
There is an upper bound on the indices of the coefficients $\cnmN$ of
eq.~(\ref{2.11}):
\qq
\intg{m_j \over 2} \leq n - \sum_{i=1}^N \intg{m_i \over 2},
\label{i2.1}
\qqq
here $\intg{m\over 2}$ denotes the integer part of ${m\over 2}$.
\end{proposition}
The proof is completely similar to that of the Proposition~3.4
of~\cite{Ro4}. Suppose that there is a coefficient $\cnmN$ for
which~(\ref{i2.1}) is not true, say, for $m_1$. If all indices $m_j$
of $\cnmN$ are odd, then the highest degree monomial of the
corresponding polynomial
\qq
\cnmN \polN
\label{2.p1}
\qqq
violates the
inequality~(\ref{2.10}). Therefore it has to be canceled by monomials
of other polynomials
\qq
\cnmNp \polNp
\label{2.p2}
\qqq
for which
\qq
m\p_j \geq m_j,\;\;
\ojN,
\qquad \sjN m_j\p > \sjN m_j.
\label{ii2.1}
\qqq
If some $m_j$ are even, then the
highest degree monomial of the polynomial~(\ref{2.p1})
is incompatible with the
structure of the \lhs of eq.~(\ref{2.11}) and it also has to be
canceled. The inequalities~(\ref{ii2.1}) show that the index $m_1$ of
the polynomials~(\ref{2.p1}) again violates~(\ref{i2.1}), so we
need to go to higher values of $\sjN m_j$ for further cancelation.
Since $\sjN m_j\leq n$, this process can not be completed. This
contradiction proves the proposition.

Now we begin to prove the Theorem~\ref{t1.1}. Following~\cite{Mu1},
we use the relations
\qq
\sak \spq^{\a^2} & = &
e^{ i{\pi \over 4} (1 - \kappa)} \sqrt{K},
\label{2.13}
\\
\sak \spq^{\a^2} & = &
\xx^{\kt} u^{-1}, \qquad
u,\,u^{-1} \in \zq
\label{2.14}
\qqq
in order to present $\knt$ in the following form:
\qq
\knt = e^{ i {\pi\over 4} N (1 - \kappa)}
\frac {u^N} {\xx^{ N \kt}}.
\label{2.15}
\qqq
Substituting this expression into eq.~(\ref{2.7}) we find that
\qq
\zpmk & = & e^{ {i\pi\over 4} (\kappa - 1)
\sjN (\sign{\xp_j} - 1) }
\left( \pjN \sign{\xp_j} \right)
\spq^{ 4^* \sjN ( 3\sign{\xp_j} - \xp_j - \xp_j^*) }
\nonumber\\
&&\times
{u^N \over \xx^{N\kt}}
\spajk \spq^{ 4^* \sjN \xp_j \a_j^2 }
\Japslk.
\label{2.16}
\qqq
Since
${1\over 4} (\kappa - 1) (\sign{\xp_j} - 1) \in \ZZ$,
we conclude that to prove Theorem~\ref{t1.1} it is enough to show
that
\qq
\spajk \spq^{4^* \sjN \xp_j \a_j^2}
\Japslk = \xx^{N\kt} w, \qquad w \in \zq.
\label{2.17}
\qqq
We will substitute the expansion~(\ref{2.8}) and check the
property~(\ref{2.17}) for every polynomial $\Dtc$ separately. The
remainder term
$\xx^{ (N + 1) \kt} \Jalkr$
obviously satisfies eq.~(\ref{2.17}). Moreover, for some
$w\in\zq$,
\qq
\left[ {1\over \xx^{N \kt} } \spajk \spq^{4^* \sjN \xp_j \a_j^2}
\Jalkr\, \xx^{ (N + 1)\kt + 1 }
\right]^{\diamond}
= \left[ \xx^{\kpt} w \right]^{\diamond}=0,
\label{2.18}
\qqq
so that we can neglect the contribution of this term in all further
calculations.

To estimate the contribution of a polynomial
$\Dtc \, \xx^n$ we need the
following simple lemma:
\begin{lemma}
\label{l2.1}
For $\xp,m \in \ZZ$, $m\geq 0$
\qq
\spak \spq^{\xp\a^2} \a^{2m} =
\xx^{\mmax{0, \kt - m}} w,\qquad w\in\zq.
\label{2.19}
\qqq
\end{lemma}
The lemma needs a proof only for $m < \kt$. To prove that an element
$u \in \zq$ is divisible by $\xx^n$, $n \leq K-1$ one may present it
as an integer coefficient polynomial of $\xx$ and check that the
coefficients in front of all $\xx^{n\p}$, $n\p < n$ are divisible by
$K$. We substitute $\spq = \xx + 1$ in eq.~(\ref{2.19}) and express
the powers of $\spq$ in terms of ``binomial'' polynomials
\qq
\!\!\!\!\!\!
\spak \spq^{ \xp \a^2} \a^{2m} =
\sum_{n\geq 0} \spak P_n(\xp\a^2) \a^{2m} \xx^n =
\sum_{n\geq 0} \sum_{0 \leq m\p \leq n}
\spak
{C_{n,m\p} \over n!}
\xp^{2m\p}
\a^{2(m + m\p)} \xx^n,
\label{2.20}
\qqq
here
\qq
\a (\a - 1) \cdots
(\a - n + 1) = \sum_{m\p \leq n} C_{n,m\p} \a^{m\p},
\qquad C_{n,m\p} \in \ZZ.
\label{2.020}
\qqq
It is well known in number theory that
\qq
\spak \a^{2m} = 0 \mod{K} \qquad \mbox{for}\;\;\;
0 \leq m < K - 2.
\label{2.21}
\qqq
Therefore
the numerator of
the contribution of a term
${C_{n,m\p} \over n!} \xp^{2m\p}\a^{2(m + m\p)} \xx^n$
will be divisible by $K$ for all
$m + m\p < \kt \leq K - 2$,
that is, for all $n < \kt - m$. The denominator $n!$ is harmless,
because since $n < K$, it is not divisible by $K$ and can not cancel
the factors of $K$ coming from the sum over $\a$. This proves the
lemma.

This lemma can be easily generalized to the ``binomial''
polynomials~(\ref{2.12}):
\begin{lemma}
\label{l2.2}
For $\xp,\xp\p,m \in \ZZ$, $\xp\p \in 2\ZZ$, $m \geq 0$
\qq
\spak \spq^{\xp\a^2}
P_m \left( {\a + \xp\p - 1 \over 2} \right) =
\xx^{ \mmax{0, \kt - \intg{m\over 2} } }\, w,
\qquad w \in \zq.
\label{2.22}
\qqq
\end{lemma}
The proof is similar to that of the previous lemma. The choice of
summation range for $\a$ obviates the fact that odd powers of $\a$ in
$P_m \left( {\a + \xp\p - 1 \over 2} \right)$
can be ignored. The numerators of the
contributions of even powers of $\a$ are
divisible by $K$. The coefficients of $P_m$ have a denominator
$m!$, but we may assume that $m < K - 1$ (otherwise eq.~(\ref{2.22})
is obvious) so that the denominator does not cancel the factors
of $K$.

Let $\mb(n)$ be the maximum value of $m$ appearing in the \lhs of
eq.~(\ref{2.11}). Then for every coefficient
$\cnmN$ from the \rhs of
that equation
\qq
m_1 + \cdots + m_N \leq \mb(n).
\label{2.23}
\qqq
Therefore we can combine eqs.~(\ref{2.11}) and~(\ref{2.22}) into the
following estimate of the contribution of polynomials $\Dtmn$ to the
sum~(\ref{2.17}):
\qq
\lefteqn{
\spak \spq^{4^* \sjN \xp_j \a_j^2} \Dtcps\, \xx^n
}
\hspace*{3in}
\nonumber\\
&&
=
\xx^{N\kt + n - \mb(n)}\, w,
\qquad w \in \zq.
\label{2.24}
\qqq
Since
$\mb(n) \leq {3\over 4} n \leq n$,
this estimate is enough to prove eq.~(\ref{2.17}) and also the
Theorem~\ref{t1.1}.
Note that the proof required only a weaker bound $m\leq n$
for $\Dtmn$ rather than a stronger bound $m\leq {3\over 4}n$
of~\cite{Ro4}. However the bound $m\leq {3\over 4}n$ is necessary to
prove that only a finite number of polynomials $\Dtmn$ contribute to
the coefficients of $\xx^{n\p}$, $n\p \leq \kt$ in the expansion of
$[\zpmk]^{\diamond}$. Indeed, since $\mb(n) \leq {3\over 4} n$, then
$n - \mb(n) \geq {1\over 4} n$ and eq.~(\ref{2.24}) suggests that we
may limit our attention to only those polynomials~(\ref{2.p1}) for
which
\qq
n \leq 2 (K-1).
\label{2.25}
\qqq
%

\nsection{Gaussian Sums and Integrals}
\label{s3}
We are going to derive a surgery formula for $[\zpmk]^{\diamond}$
which would express it in terms of the derivatives $\Dtmn$ of the
colored Jones polynomial. As we will see, this requires a calculation
of the gaussian sum
\qq
\gpq{m} = \eikap
\spak \spq^{\xq^* p \a^2} \a^{2m} \xx^m,
\qquad m\leq \kt.
\label{3.2}
\qqq
More precisely, we need to find only $[\gpq{m}]^{\diamond}$. We
already know that $\gpq{m} \in \zq$.
\begin{proposition}
\label{p3.1}
The sum of eq.~(\ref{3.2}) is related to the gaussian integral. For
$\xp,\xq \in \ZZ$, $\xp,\xq \neq 0 \mod{K}$, $0\leq m \leq \kt$,
\qq
[\gpq{m}]^\diamond & = &
\Leg{\xp\xq^*} \left[ e^{ -i{\pi\over 4} \sign{\xp\over \xq}}
\left( {2\over K} \right)^{1\over 2}
\left| {\xp\over \xq} \right|^{1\over 2} \int_{-\infty}^{+\infty}
d\a\, e^{ {2\pi i\over K} {\xp\over \xq} \a^2}
\a^{2m}\,\xx^m \right]^\vee
\label{3.3}\\
&& \qquad + \xx^{\kpt - m} w,\;\;w\in\zq.
\nonumber
\qqq
\end{proposition}

To prove the proposition we calculate the following sum:
\qq
\eikap \spak \spq^{ \xp\xq^* \a^2} \spq^{2n\a} =
\eikap \spq^{-\xp^* q n^2} \spak
\spq^{\xp \xq^* (\a + n \xp^* \xq)^2 }.
\label{3.4}
\qqq
Since (for $\xp^* \in \ZZ$)
\qq
\spak
\spq^{\xp \xq^* (\a + n \xp^* \xq)^2 }
& = & \spak \spq^{\xp \xq^* \a^2}
= \sta \spq^{\xp \xq^* (2\ta + 1)^2}
\nonumber
\\
& = & \sta \spq^{4\xp\xq^* (\ta + 2^*)^2}
= \sqrt{K} e^{ i{\pi\over 4} (1 - \kappa)} \Leg{\xp\xq^*},
\label{3.5}
\qqq
we find that
\qq
\eikap \spak \spq^{ \xp\xq^* \a^2} \spq^{2n\a} =
\Leg{\xp\xq^*} \spq^{-\xp^* \xq n^2}.
\label{3.6}
\qqq
We substitute $\spq = 1 + \xx$ in $\spq^{2n\a}$ and $\spq^{-\xp^*\xq
n^2}$. After going from $\zq$ to the factor-ring
$\zpkx = \zkx/ \xx^{\kpt} \zkx$ and using the ``checked binomial
polynomial''
\qq
P_m^\vee(\a) = (m!)^* \a (\a - 1) \cdots (\a - m + 1)
= (m!)^*\sum_{l=0}^m C_{m,l} \a^l,
\label{3.7}
\qqq
we find that
\qq
\left[ \sum_{m=0}^{K-1} \sum_{l=0}^{\intg{m\over 2}}
(m!)^* C_{m,2l} \gpq{l} 4^l n^{2l} \xx^{m-l} \right]^\diamond
=
\Leg{\xp\xq^*} \left[ \sum_{m=0}^{\kt}
\xx^m
\sum_{l=0}^{m}
(m!)^* C_{m,l} (-1)^l (\xp^* \xq)^l n^{2l} \right]^\diamond
\label{3.8}
\qqq
We limited the sum over $m$ in the \lhs of this equation to $m\leq
K-1$ because for $m\geq K$ the minimum power of $\xx^{m-l}$ is
greater that $\kt$. Note that since $m\leq K-1$, then $(m!)^*$ is
well defined.

If we substitute the expansion
\qq
[\gpq{l}]^\diamond =
\sum_{m=0}^{\kt - l} \gpqi{l}{m}\,\xx^m +
\xx^{\kpt - l}w,\qquad w\in\zq
\label{3.9}
\qqq
into eq.~(\ref{3.8}), then we can find all the coefficients
$\gpqi{l}{m}$ by equating the coefficients of \lhs and \rhs of
eq.~(\ref{3.8}) at equal powers of $\xx$ and $n$. These coefficients
have to be equal due to the following simple lemma:
\begin{lemma}
\label{l3.1}
If a degree of a polynomial $P(n) \in \ZZ_K[n]$ is less than $K$ and
$P(n)=0$ for all $n\in \ZZ_K$, then all the coefficients of $P(n)$
are zero modulo $K$.
\end{lemma}

The proof follows from the fact that the $K\times K$ Van-der-Monde
determinant in $\ZZ_K$ is non-zero.

Each gaussian sum $\gpq{l}$ appears in the \lhs of eq.~(\ref{3.8})
with its own power of $n$: $n^{2l}$. Therefore the coefficients
$\gpqi{l}{m\p}$ of eq.~(\ref{3.9}) can be calculated by ``dividing''
the polynomial
$(-1)^l (\xp^* \xq)^l \sum_{m=l}^{\kt} (m!)^* C_{m,l} \xx^m$
appearing at $n^{2l}$ in the \rhs of eq.~(\ref{3.8}) by the
polynomial
$\sum_{m=2l}^{K-1} C_{m,2l} 4^l \xx^{m-l}$
appearing in the \rhs of that equation at the same power on $n$. The
division is not quite well-defined, hence the indeterminacy in the
elements $w$ of eq.~(\ref{3.9}).

This whole calculation of dividing the polynomials can be made more
explicit if we go back to eq.~(\ref{3.6}) and make the following
substitutions:
\qq
\spq^{2n\a} & = &
e^{2n\a \log (1 + \xx)} = \sum_{l=0}^\infty
{(2n\a)^l \over l!} (\log(1 + \xx))^l,
\label{3.10}
\\
\spq^{-\xp\xq^* n^2} & = &
e^{-\xp\xq^* n^2 \log(1 + \xx)} = \sum_{l=0}^{\infty}
(-1)^l {(\xp\xq^* n^2)^l \over l!} (\log(1 + \xx))^l.
\label{3.11}
\qqq
After ``checking'' the logarithm
\qq
\log^\vee (1 + \xx) = \xx \sum_{n=0}^{K-2} (-1)^n (n+1)^* \xx^n
\label{3.12}
\qqq
we see that eq.~(\ref{3.8}) transforms into
\qq
\lefteqn{
\left[ \sum_{l=0}^{\kt} \gpq{l} (2n)^{2l} ((2l)!)^*
(\log^\vee(1 + \xx))^l
\left( { \log^\vee(1 + \xx)\over \xx} \right)^l
\right]^\diamond
}\hspace*{2in}
\label{3.13}
\\
&&
=
\Leg{\xp\xq^*} \left[ \sum_{l=0}^{\kt} (-1)^l n^{2l}
(\xp\xq^*)^l (l!)^* (\log^\vee(1 + \xx))^l \right]^\diamond.
\nonumber
\qqq
Thus we find that
\qq
[\gpq{m}]^\diamond & = &
(-1)^m \Leg{\xp\xq^*} (2^*)^{2m} (2m)! (m!)^*
\left[ \left( {\xx\over \log^\vee (1+\xx)} \right)^m \right]^\diamond
\nonumber\\
&&\qquad\qquad
+ \xx^{\kpt - m} w,
\qquad w\in \zq.
\label{3.14}
\qqq

Consider now the following identity which is an integral
analog of eq.~(\ref{3.6}):
\qq
\int_{-\infty}^{+\infty} d\a\, \spq^{ {\xp\over \xq} \a^2}
\spq^{2n\a} =
e^{i{\pi\over 4}\sign{\xp\over \xq}}
\left({K\over 2}\right)^{1\over 2}
\left| {\xq\over \xp} \right|^{1\over 2}
\spq^{ -{\xq\over \xp} n^2}.
\label{3.15}
\qqq
After substituting eqs.~(\ref{3.10}),~(\ref{3.11}) we find that
\qq
\int_{-\infty}^{+\infty} d\a\, \spq^{ {\xp\over \xq} \a^2}
\a^{2m} \xx^m
=
e^{i{\pi\over 4}\sign{\xp\over \xq}}
\left({K\over 2}\right)^{1\over 2}
\left| {\xq\over \xp} \right|^{1\over 2}
{ (-1)^m (2m)! \over 4^m m!}
\left( {\xx \over \log(1 + \xx)} \right)^m
\label{3.16}
\qqq
Eq.~(\ref{3.3}) follows from comparing eq.~(\ref{3.16}) to
eq.~(\ref{3.14}).

The formula~(\ref{3.3}) can be generalized to the type of summands
that appear in eq.~(\ref{2.7}) after the substitutions~(\ref{2.8})
and~(\ref{2.11}).
\begin{proposition}
\label{p3.2}
For $\xp \neq 0  \mod{K}$, $0\leq m\leq K$
\qq
\lefteqn{
\left[ \eikap \spak \spq^{4^* \xp \a^2}
\pmaps \xx^{\intg{m\over 2}}
\right]^\diamond
}
\label{3.17}
\\
= &&\Leg{\xp}
\left[ e^{ -i{\pi\over 4} \sign{\xq}}
\left( {|p|\over 2K} \right) ^{1\over 2}
\int_{-\infty}^{+\infty} d\a\,
e^{{i\pi\over 2K} \xp \a^2}
\pmap \xx^{\intg{m\over 2}}
\right]^\diamond
+ \xx^{\kpt - \intg{m\over 2}} w,
\nonumber\\
&&\qquad w\in \zq.
\nonumber
\qqq
\end{proposition}

To prove the proposition for $m < K$ we substitute
$P_m^\vee(2^*(\a - 1 + \xp^*))$
for $\pmaps$ in the \lhs of this
equation and then take the sum for each monomial of
$P_m^\vee(2^*(\a - 1 + \xp^*))$
separately. If $m$ is odd then the highest power $\a^m$ does not
contribute to the sum, therefore the factor
$\xx^{\intg{m\over 2}}$
is enough to apply eq.~(\ref{3.3}). We also used the multiplicativity
of Legendre symbol and $\Leg{4^*}=1$, since $4^* = (2^*)^2$.

The case of $m=K$ requires a special care. We start with the \lhs of
eq.~(\ref{3.17}). We can use the symmetry of the summation range and
gaussian exponent in order to substitute
\qq
\pkev(\a) = {1\over 2} \left( \pKaps + \pKapsm \right)
\label{3.18}
\qqq
instead of $\pKaps$. The even polynomial~(\ref{3.18}) takes integer
values for odd $\a$ and its degree is equal to $K-1$. Therefore the
highest divisor of denominators of its coefficients is $K-1$ and we
can apply all our previous results to the calculation of the
contribution of its monomials. Eq.~(\ref{3.17}) indicates that we
need to determine only the terms of order $\xx^0$, hence we are
interested only in the contribution of the highest degree monomial
\qq
{\xp^* - K\over 2^K \,(K-1)!} \a^{\kt},
\label{3.1018}
\qqq
which is determined with the help of eq.~(\ref{3.14}). Consider now
the \rhs of eq.~(\ref{3.17}). Again we substitute the even polynomial
\qq
\pkev(\a) = {1\over 2} \left( \pKap + \pKapm \right).
\label{3.018}
\qqq
Some of its denominators may have $K$ as a divisor, but according to
eq.~(\ref{3.17}) we are interested only in the contribution of the
highest power of $\a$:
\qq
{ {1\over p} - K \over 2^K\, (K-1)!} \a^{\kt}.
\label{3.2018}
\qqq
Comparing it to monomial~(\ref{3.1018}) and applying eq.~(\ref{3.3})
to their contributions we arrive at eq.~(\ref{3.17}). This ends the
proof of the Proposition~\ref{p3.2}.

Next we move to the polynomials $\Dtmn$ which participate in the
expansion eq.~(\ref{2.8}):
\begin{proposition}
\label{p3.3}
The gaussian sums and integrals of polynomials $\Dtmn$ are related by
the equation
\qq
\lefteqn{
\left[ \frac
{e^{ i {\pi\over 4} N (\kappa - 1)}} {K^{N\over 2}}
\spajk \spq^{4^* \sjN \xp_j \a_j^2}
\Dtcps\, \xx^n
\right]^\diamond
}
\hspace*{1in}
\nonumber\\
&
=&
 \Leg{\pjN \xp_j}
\left[
e^{-i{\pi\over 4} \sjN \sign{\xp_j}} (2K)^{-{N\over 2}}
\left|\pjN \xp_j \right|^{1\over 2}
\right.
\label{3.19}\\
&&
\hspace*{-1in}
\times
\left.
\int_{-\infty}^{+\infty} d\a\,
e^{ {i\pi\over 2K} \sjN \xp_j \a_j^2} \Dtcp
\right]^\diamond.
\nonumber
\qqq
\end{proposition}

To prove this proposition we rearrange the
representation~(\ref{2.11}) in the following form:
\qq
\lefteqn{
\Dtc\,\xx^n
}
\label{3.20}
\\
&&
=
\xx^{n-\sjN \intg{m_j\over 2} }
\smb \cnmN \xx^{\intg{m_1\over 2}}
P_{m_1}\left( {\a_1 - 1\over 2} \right)
\cdots
\xx^{\intg{m_N\over 2}}
P_{m_N}\left( {\a_N - 1\over 2} \right).
\nonumber
\qqq
We know from Lemma~\ref{l2.2} that the contribution of each
polynomial
$\xx^{\intg{m_j\over 2}} P_{m_j}
\left( { \a_j - 1\over 2}\right)$
to the \lhs of eq.~(\ref{3.19}) belongs to $\zq$. Therefore
the contribution of the whole expression~(\ref{3.20}) starts at
$\xx^{n - \sjN \intg{m_j\over 2}}$. Hence we may assume that
\qq
n - \sjN \intg{m_j\over 2} \leq \kt,
\label{3.21}
\qqq
otherwise the contribution of the polynomial
\qq
\xx^{n-\sjN \intg{m_j\over 2} }
\cnmN \xx^{\intg{m_1\over 2}}
P_{m_1}\left( {\a_1 - 1\over 2} \right)
\cdots
\xx^{\intg{m_N\over 2}}
P_{m_N}\left( {\a_N - 1\over 2} \right).
\label{3.22}
\qqq
is annihilated by the homomorphism $^\diamond$.

The inequalities~(\ref{i2.1}) and~(\ref{3.21}) mean that
$m_j\leq\kt$, so we can apply the Proposition~\ref{p3.2} to every
polynomial
$\xx^{\intg{m_j\over 2}} P_{m_j}
\left( { \a_j - 1\over 2}\right)$.
The terms
$
\xx^{\kpt - \intg{m_j\over 2}}w$ of eq.~(\ref{3.17}) can be neglected
because
\qq
\left[ \xx^{n - \sum_{i=1}^N \intg{m_i\over 2}}\,
\xx^{\kpt - \intg{m_j\over 2}}
\,w
\right]^\diamond
=
\left[
\xx^{\kpt + \left(
n - \sum_{i=1}^N \intg{m_i\over 2} - \intg{m_j\over 2}
\right) }
\,w
\right]^\diamond
=0
\label{3.23}
\qqq
in view of the inequality~(\ref{i2.1}). This proves the
Proposition~\ref{3.3}.

Now we can prove the Proposition~\ref{p1.1}. We substitute
eq.~(\ref{2.8}) into eq.~(\ref{2.4}), apply the homomorphism
$^\diamond$ and retain only the relevant terms from the sum of
eq.~(\ref{2.8}). Since the contribution of the \lhs of
eq.~(\ref{3.20}) starts at
$\xx^{n - \sjN\intg{m_j\over 2}}$ and
$\sjN \intg{m_j\over 2} \leq m \leq {3\over 4}n$, it is enough to
retain only the terms with
$n\leq 2(K-1)$:
\qq
[\zpmk]^\diamond
& = &
\left[
e^{i {\pi\over 4} (\kappa - 1) \sjN (\sign{\xp_j} - 1) }
\left( \pjN \sign{\xp_j} \right)
\spq^{ 4^* \sjN (3\sign{\xp_j} - \xp_j - \xp_j^*) }
\,
\frac { e^{i {\pi\over 4} N (\kappa - 1)}} {K^{N\over 2}}
\right.
\nonumber
\\
&&
\hspace*{-1in}
\left.
\times
\spajk \spq^{4^* \sjN \xp_j \a_j^2}
\ppajps
\sum_{n=0}^{2(K-1)}
\smnf \Dtmn (\aaps) \xx^n
\right]^\diamond
\nonumber
\\
& = &
\Leg{\pjN \xp_j}
e^{i {\pi\over 4} (\kappa - 1) \sjN (\sign{\xp_j} - 1) }
\label{3.24}
\\
&&\times
\left[
e^{-i{\pi\over 4}\sjN \sign{\xp_j}}
\left( \pjN \sign{\xp_j} \right)
\spq^{ {1\over 4} \sjN (3\sign{\xp_j} - \xp_j - \xp_j^*) }
\frac {\left| \pjN \xp_j \right|^{1\over 2} }
{(2K)^{N\over 2}}
\right.
\nonumber
\\
&&
\times
\int_{-\infty}^{+\infty} d\a_1\cdots d\a_N\,
e^{{i\pi\over 2K} \sjN \xp_j \a_j^2}
\nonumber\\
&&
\qquad
\times\left.
\ppajp
\sum_{n=0}^{2(K-1)}
\smnf \Dtmn (\aap) \xx^n
\right]^\vee.
\nonumber
\qqq
Here we used an identity
\qq
\left[
\spq^{ 4^* \sjN (3\sign{\xp_j} - \xp_j - \xp_j^*) }
\right]^\diamond
=
\left[
\spq^{ {1\over 4} \sjN (3\sign{\xp_j} - \xp_j - \xp_j^*) }
\right]^\vee.
\label{3.25}
\qqq

We can extend the sum over $n$ in the \rhs of eq.~(\ref{3.25}) to all
$n\geq 0$ because the contribution of $\Dtmn \xx^n$ with $n>2(K-1)$
starts above $\xx^{\kt}$. As a result, we obtain the full Jones
polynomial. Now using the identities
\qq
&\displaystyle
e^{i {\pi\over 4} (\kappa - 1) (\sign{\xp} - 1)}
\Leg{\xp} = \Leg{|p|},
\label{3.26}\\
&\displaystyle
\left|\pjN \xp_j \right| = \ordH
\label{3.27}\\
&\displaystyle
\sign{\xp} = i e^{- i{\pi\over 2}\sign{\xp}},
\label{3.28}
\qqq
we get the following formula
\qq
[\zpmk]^\diamond
& = &
\Leg{\ordH}
\left[
e^{-{3\over 4} i\pi \sjN \sign{\xp_j} }
{i^N\over (2K)^{N\over 2}}
e^{ {i\pi\over 2K} \sjN
\left(
\sign{\xp_j} - \xp_j - {1\over \xp_j}
\right) }
\left|\pjN \xp_j \right|^{1\over 2}
\right.
\nonumber
\\
&&\left.
\qquad
\times
\spint{\a_j} d\a_1\cdots d\a_N \,
e^{ {i\pi\over 2K} \sjN \xp_j\a_j^2}
\Japlk
\right]^\vee
\label{3.29}
\qqq
Combining it with eqs.~(\ref{1.25}) and~(\ref{1.17}) and using
eq.~(\ref{3.27}) again we easily arrive at eq.~(\ref{1.21}). Note
that eq.~(\ref{3.27}) guarantees that $\xp_j\neq 0\mod{K}$ if
$\ordH \neq 0 \mod{K}$. This proves the Proposition~\ref{p1.1}.

\nsection{A General Rational Surgery Formula}
\label{s4}
Up until this point we were working only with surgeries of the type
$(\xp,1)$. This was enough to prove the theorems of Section~\ref{s1},
because H.~Murakami and T.~Ohtsuki showed~\cite{Mu2} that any \Rhs
$M$ can be constructed by $(\xp_j,1)$ surgeries on an \asl in $S^3$
up to a connected sum of lens spaces $L_{\xp\p_j,1}$, that is,
instead of $M$ one might end up with
$M\# L_{\xp\p_1,1} \# \ldots L_{\xp\p_n,1}$. However from the
technical point of view it would be better to have a formula for
the invariant $\zpmk$ of a manifold constructed by general rational
surgeries $(\xp_j,\xq_j)$ on the components of an $N$-component link
$\cL$ in $S^3$.

The formula for Witten's invariant $\zmk$ was derived by
L.~Jeffrey~\cite{Je}
\qq
\zmk & = &
\zsk \exp \left[
i {\pi\over 4} {K - 2\over K}
\left( \sjN \Phi(U^{(\xp_j,\xq_j)} - 3 \sign{L} \right) \right]
\label{4.1}
\\
&&\qquad
\times
\sum_{1\leq\a_1,\ldots,\a_N \leq K-1}
\Jalk \pjN \tU^{(\xp_j,\xq_j)}_{\a_j 1},
\nonumber
\qqq
here $L$ is the linking matrix of $\cL$, ${\xp_j\over \xq_j}$ being
the self-linking numbers. The matrices
\qq
U^{(\xp_j,\xq_j)} =
\pmatrix { \xp_j & \xr_j \cr \xq_j & \xs_j\cr} \in SL(2,\ZZ)
\label{4.2}
\qqq
describe the surgeries (a meridian on the tubular neighborhood is
glued to $\xp_j\mbox{(meridian)} + \xq_j \mbox{(parallel)}$ of the
link complement),
\qq
\Phi
\left[
\matrix
{ \xp & \xr \cr \xq & \xs\cr}
\right] =
{\xp + \xs \over \xq} - 12 s(\xp,\xq),
\label{4.3}
\qqq
$s(\xp,\xq)$ being the Dedekind sum, and
\qq
\tU^{(\xp,\xq)}_{\a\b} & = &
i {\sign{\xq}\over \sqrt{2K|q|}}
e^{ -i{\pi\over 4} \Phi (U^{(\xp,\xq)})}
\sum_{n=0}^{q-1} \sum_{\mu = \pm 1} \mu
\nonumber\\
&&\qquad\times
\exp \left[ {i\pi\over 2K\xq}
(\xp\a^2 - 2\a(2Kn + \mu\b) + \xs (2Kn + \mu \b)^2) \right].
\label{4.4}
\qqq

Let us introduce some notations. A rational $(\xp_j,\xq_j)$ surgery
on $\cL_j$ can be presented as a combination of $(m_t^{(j)},1)$,
$1\leq t \leq \tb^{(j)}$ surgeries on a chain of unknots simply
linked to $\cL_j$ (see \eg~\cite{FrGo},~\cite{Je} and references
therein) such that
\qq
U^{(\xp_j,\xq_j)} =
T^{m^{(j)}_{\tb^{(j)}}} S
T^{m^{(j)}_{\tb^{(j)}-1}} S
\cdots
T^{m_1^{(j)}} S,
\label{4.5}
\qqq
here
\qq
S=\pmatrix{ 0&-1\cr 1&0\cr},
\qquad T=\pmatrix{ 1& 1\cr 0&1\cr}.
\label{4.05}
\qqq
We denote this chain (including $\cL_j$ itself) as $\tcL_j$ and all
the chains of $\cL$ as $\tcL$. For $1\leq t \leq \tb^{(j)}$ we set
\qq
U^{(\xp_t^{(j)},\xq_t^{(j)})} =
T^{m_t^{(j)}} S \cdots T^{m_1^{(j)}} S,
\label{4.6}
\qqq
so that $\xp_{\tb^{(j)}}^{(j)} = \xp_j$,
$\xq_{\tb^{(j)}}^{(j)} = \xq_j$. From now on we assume for simplicity
that none of the numbers $\xq_t^{(j)}$ is divisible by $K$. Then we
are going to prove the following:
\begin{proposition}
\label{p4.1}
Let $M$ be a manifold constructed by $(\xp_j,\xq_j)$ surgeries on an
$N$-component link $\cL$ in $S^3$. Then
\qq
\zpmk & = &
\Leg{ \left| \pjN \xq_j \right| }
\frac {\pjN \sign{\xq_j} } {K^{N\over 2}}
e^{-i {\pi\over 4} \kappa \sign{L} }
e^{- i\pi {3\over 4} {K-2\over K} \sign{L} }
\nonumber\\
&&
\qquad
\times
\spq^{-4^* \sjN \Phi ( U^{(\xp_j,\xq_j)} ) }
\spq^{\left( 2^* - {1\over 2} \right) \sjN \sign{\xp_j\over \xq_j}}
\label{4.7}
\\
&&
\qquad
\times
\spajk \Jalk\,
\spq^{4^*\sjN \xq_j^* (\xp_j\a_j^2 + \xs_j) }
\pjN \left( {i\over 2} \smj \mu_j
\spq^{ -2^* \xq_j^* \mu_j\a_j }
\right)
\nonumber
\qqq
\end{proposition}

We could use the general surgery formula~(\ref{4.7}) instead of the
$(\xp_j,1)$ surgery formula~(\ref{1.5}) throughout the
Sections~\ref{s2} and~\ref{s3} in order to produce a somewhat more
flexible proof of Theorems~\ref{t1.1},~\ref{t1.2} and
Proposition~\ref{p1.1}.

We begin the proof of the Proposition~\ref{p4.1}
by recalling the Kirby-Melvin formula~\cite{KiMe}
which expresses $\zpmk$ in terms of data associated to $\tcL$:
\qq
\zpmk & = &
e^{i {\pi\over 4} {K-2\over K} \left[
\sjN \Phi(U^{(\xp_j,\xq_j)} ) - 3 \sign{L} \right] }
\,
e^{- i {\pi\over 4} \kappa \sign{\tL} }
\nonumber
\\
&&\qquad\times
\spajk \Jalk
\pjN \chU^{(\xp_j,\xq_j)}_{\a_j 1},
\label{4.8}
\qqq
here (we drop the index $j$ in eq.~(\ref{4.5}) )
\qq
&\displaystyle
\chU^{(\xp,\xq)}_{\a\b} =
\left( \chT^{m_{\tb}}\chS\,\chT^{m_{\tb-1}} \chS
\cdots \chT^{m_1} \chS \right)_{\a\b},
\label{4.9}
\\
&\displaystyle
\chT_{\a\b} = e^{-i{\pi\over 4}}
\spq^{{1\over 4}\a^2} \delta_{\a\b},
\qquad
\chS_{\a\b} = {1\over \sqrt{K}} \sin \left(
{\pi\over K} \a\b \right)
\label{4.10}
\qqq
and $\tL$ is the linking matrix of the ``expanded'' link $\tcL$.

The following lemma presents an explicit expression for
$\chU^{(\xp,\xq)}_{\a\b}$, which is similar to that of
eq.~(\ref{4.4}):
\begin{lemma}
\label{l4.1}
\qq
\chU^{(\xp,\xp)}_{\a\b} & = &
\Leg{|\xq|} {\sign{\xq}\over \sqrt{K}}
e^{-i{\pi\over 4} \Phi(U^{(\xp,\xq)} )}
\,e^{ -i {\pi\over 4} \kappa \sum_{t=0}^{\tb-1}
\sign{\xp_t\over \xq_t} }
\label{4.11}
\\
&&\qquad\times
\spq^{ \left({1\over 4} - 4^* \right) \sum_{t=1}^{\tb} m_t
\,
+
\,
\left( 2^* - {1\over 2} \right) \sum_{t=1}^{\tb}
\sign{\xp_t\over \xq_t} }
\,
\spq^{4^* \xq^* (\xp\a^2 + \xs \b^2)}
\,
{i\over 2}\sum_{\mu = \pm 1} \mu
\,\spq^{-2^*\xq^* \mu\a\b}.
\nonumber
\qqq
\end{lemma}
To prove the lemma we slightly change eq.~(\ref{4.10}):
\qq
\chT_{\a\b} = e^{-i{\pi\over 4}} \,
\spq^{{1\over 4} - 4^*} \,
\spq^{4^* \a^2}\,\delta_{\a\b},\qquad
\chS_{\a\b} = {1\over \sqrt{K}} {i\over 2}
\,\spq^{\pm\left( 2^* - {1\over 2} \right) }
\sum_{\mu=\pm 1} \mu\,\spq^{-2^*\mu\a\b},
\label{4.12}
\qqq
the choice of sign in
$\spq^{\pm \left( 2^* - {1\over 2} \right)}$
is arbitrary.

We prove eq.~(\ref{4.11}) by induction on $\tb$. If $\tb=1$, that is,
if $U^{(\xp,\xq)} = T^m S$, then the check is trivial if we recall
that $\Phi(T^m S) = m$. It is also easy to check eq.~(\ref{4.11}) for
$U^{(\xp+m\xq,\xq)} = T^m U^{(\xp,\xq)}$.
It remains to check
eq.~(\ref{4.11}) for
$U^{(-\xq,\xp)} = S U^{(\xp,\xq)}$.

We have to calculate the sum:
\qq
\chU^{(-\xq,\xp)}_{\a\b} =
\spgk \chS_{\a\gamma} \chU^{(\xp,\xq)}_{\gamma\b}.
\label{4.13}
\qqq
The following gaussian sum is at the center of this calculation:
\qq
\lefteqn{
\spgk \sum_{\mu_{1,2}=\pm 1} \mu_1 \mu_2
\,
\spq^{4^* \xq^* \xp \gamma^2 - 2^* \gamma (\xq^* \mu_1 \b + \mu_2\a)
+ 4^* \xq^* \xs \b^2 }
}
\hspace*{1in}
\label{4.14}
\\
&=&
\spgk \sum_{\mu_{1,2}=\pm 1} \mu_1 \mu_2
\,
\spq^{4^* \xq^* \xp ( \gamma - \xp^* \mu_1\b - \xp^* \xq \mu_2 \a )^2
- 4^* \xq^* \xp^* (\mu_1\b + \xq \mu_2 \a)^2 + 4^* \xq^* \xs \b^2 }
\nonumber\\
& = &
2 \sqrt{K} \eikapw \Leg{\xp \xq^* }
\,
\spq^{4^* \xp^* (-\xq \a^2 + \xr \b^2) }
\sum_{\mu= \pm 1} \mu
\, \spq^{-2^* \xp^* \mu\a\b},
\nonumber
\qqq
here we used the following relations:
\qq
&\displaystyle
\spgk \spq^{4^* \xp \xq^* \gamma^2 }
= \eikapw \sqrt{K} \Leg{\xp\xq^*},
\label{4.15}
\\
&\displaystyle
\xs - \xp^* = \xp^* \xq \xr,
\label{4.16}
\qqq
the latter relation follows from $\xp\xs - \xq \xr = 1$. To complete
the verification of eq.~(\ref{4.11}) we recall the following
identities:
\qq
&\displaystyle
\Phi(SU^{(\xp,\xq)} ) = \Phi ( U^{(\xp,\xq)}) - 3\sign{\xp\over \xq},
\label{4.17}
\\
&\displaystyle
i\sign{\xq} = e^{i{\pi\over 2} \sign{\xp\over \xq} }
\sign{\xp},
\label{4.18}
\\
&\displaystyle
e^{i{\pi\over 4} (\kappa-1)} \Leg{|q|}\Leg{\xp\xq^*} =
e^{i{\pi\over 4} (\kappa - 1) \sign{\xp\over \xq} } \Leg{|p|}.
\label{4.19}
\qqq
This ends the proof of the lemma.

To finish the proof of Proposition~\ref{4.1} we rearrange some phase
factors of eqs.~(\ref{4.8}) and ~(\ref{4.11}). We substitute a
relation
\qq
\sign{\tL} = \sign{L} + \sjN \sum_{t=1}^{t^{(j)} - 1}
\sign{\xp_t^{(j)}\over \xq_t^{(j)} }
\label{4.20}
\qqq
into the factor
$e^{- i{\pi\over 4} \kappa \sign{\tL} }$
of eq.~(\ref{4.8}). Then we calculate the combination of phases
coming from that factor and from eq.~(\ref{4.11}) (we drop the index
$j$):
\qq
\lefteqn{
\spq^{\left( {1\over 4} - 4^* \right) \sum_{t=1}^{\tb} m_t
+ \left( 2^* - {1\over 2} \right)
\sum_{t=1}^{\tb}\sign{\xp_t\over\xq_t} }
\,
e^{-i{\pi\over 2}\kappa\sum_{t=1}^{\tb-1}\sign{\xp_t\over \xq_t} }
}
\label{4.21}
\hspace*{1in}
\\
&=&
\spq^{\left( {1\over 4} - 4^* \right) \sum_{t=1}^{\tb} m_t
+ \left( 2^* - {1\over 2} \right)
\sign{\xp\over\xq} }
\,
\spq^{\left(2^* - {1\over 2} - {\kappa K\over 4} \right)
\sum_{t=1}^{\tb-1} \sign{\xp_t\over \xq_t} }
\nonumber\\
& = &
\spq^{\left( {1\over 4} - 4^* \right)
\left[ \sum_{t=1}^{\tb} m_t - 3 \sum_{t=1}^{\tb-1}
\sign{\xp_t\over \xq_t} \right] }
\,
\spq^{\left( 2^* - {1\over 2} \right) \sign{\xp\over \xq} }
\nonumber\\
& = &
\spq^{ \left( {1\over 4} - 4^*\right) \Phi(U^{(\xp,\xq)} )}
\,
\spq^{ \left(2^* - {1\over 2} \right) \sign{\xp\over \xq} }
\nonumber
\qqq
We used here eq.~(\ref{2.6}) and the formula
\qq
\Phi(U^{(\xp,\xq)} ) = \sum_{t=1}^{\tb} m_t
- 3 \sum_{t=1}^{\tb - 1} \sign{\xp_t\over \xq_t}
\label{4.22}
\qqq
(see~\cite{Je} and references therein).

A combination of eqs.~(\ref{4.8}),~(\ref{4.11}),~(\ref{4.20})
and~(\ref{4.21}) leads to eq.~(\ref{4.7}). This concludes the proof
of Proposition~\ref{p4.1}.

As an application of the Proposition~\ref{p4.1} let us calculate the
invariant of a lens space $L_{\xp,\xq}$. This manifold is constructed
by a $U^{(-\xp,\xq)}$ surgery on an unknot in $S^3$. Since
\qq
J_\a(\mbox{unknot};k) = {\sin \left( {\pi\over K} \a \right)
\over \sin \left( {\pi\over K} \right) },
\label{4.23}
\qqq
we can say that $L_{\xp,\xq}$ is constructed by a chain surgery
$U^{(-\xq,-\xp)} = S U^{(-\xp,\xq)}$ applied to an empty knot, times
a factor
${ \sqrt{K}\over \sin \left( {\pi\over K} \right)}$.
Then we can read the result directly from eq.~(\ref{4.7}) by setting
there $N=1$, $\sign{L} = 0$, $\a = 1$ and
\qq
U^{(-\xq,-\xp)} =
\pmatrix{ -\xq & \xs \cr -\xp & \xr \cr}
\label{4.24}
\qqq
instead of $U^{(\xp,\xq)}$:
\qq
Z\p(L_{\xp,\xq}) & = &
\Leg{|\xp|} \sign{\xp}\,
\spq^{4^* ( \xp^* (\xq - \xr) - \Phi(U^{(-\xq,-\xp)}))}
\,
\spq^{2^*}\,
{\spq^{2^*\xp^*} - \spq^{-2^*\xp^*} \over \spq - 1}
\nonumber\\
& = &
\Leg{|\xp|} \sign{\xp} \,
\spq^{3 s^\vee (\xq,\xp)} \,
\spq^{2^*}
\,
{\spq^{2^*\xp^*} - \spq^{-2^*\xp^*} \over \spq - 1}
\label{4.25}\\
& = &
\Leg{|\xp|} \sign{\xp}
\,
\spq^{3s^\vee(\xq,\xp)}
\,
(-1)^{\xp^* + 1}
\frac {\spq^{\xp^*\over 2} - \spq^{-{\xp^*\over 2}}}
{\spq^{1\over 2} - \spq^{-{1\over 2}}},
\nonumber
\qqq
here $s^\vee(\xq,\xp)$ is the ``checked'' Dedekind sum, that is, its
denominator is inverted modulo $K$ as in eq.~(\ref{1.12}).
We used the relation
$\spq^{2^*} = - \spq^{1\over 2}$
in order to derive the last expression in this equation. Although it
might look simpler than the previous one, it obscures the fact that
$Z\p(L_{\xp,\xq}) \in \zq$. Also note that the expression
$\spq^{\xp^*\over 2}$ by itself is ambiguous since $\xp^*$ is defined
only modulo $K$.

Comparing the second expression in the \rhs of eq.~(\ref{4.25}) with
the formula for the trivial connection contribution to Witten's
invariant of the lens space
\qq
\ztr(L_{\xp,\xq};k) =
\zsk {\sign{\xp}\over \sqrt{|\xp|} } \spq^{ 3s(\xp,\xq) }
\frac {\spq^{p\over 2} - \spq^{-{p\over 2}} }
{\spq^{1\over2} - \spq^{-{1\over 2}} }
\label{i4.1}
\qqq
derived in~\cite{Je} we can easily check the relation~(\ref{1.18}).

The formula~(\ref{4.7}) also allows us to check the
Proposition~\ref{p1.1} for Seifert manifolds which are rational
homology spheres. Consider an $(N+1)$-component link $\cL$ in $S^3$
consisting of $N$ unknots $\cL_j$, $\ojN$ simply linked to a single
unknot $\cL_0$. A Seifert manifold $\seif$ is constructed by
performing $(p_j,q_j)$ surgeries on the components $\cL_j$ and a
$(0,1)$ surgery on $\cL_0$ (see, \eg~\cite{FrGo}). The Jones
polynomial of $\cL$ is known~\cite{Wi} to be equal to
\qq
J_{\b,\aa} (\cL;k) = {1\over \spk} {\pjN \spbaj \over \sin^{N-1}
\left( \pk \b \right) },
\label{i4.2}
\qqq
here $\b$ is the color of $\cL_0$ and $\a_j$ are the colors of
$\cL_j$. The signature of $\cL$ is equal to
\qq
\sign{L} = -\sph + \sjN \spqj
\label{i4.3}
\qqq
Here we introduced notations
\qq
P = \pjN p_j, \qquad
H = P \sjN {q_j\over p_j},
\label{i4.4}
\qqq
so that
\qq
\left| H_1 \left( \seif \right) \right| = |H|.
\label{i4.5}
\qqq

The Proposition~\ref{p4.1} provides the following expression for the
invariant of the Seifert manifold $\seif$:
\qq
Z\p(X;k) &=&
{i\over 2\sqrt{K} }
e^{i {\pi\over 4} \kappa \sph }
e^{i\pi {3\over 4} {K-2\over K} \sph }
\spq^{\left( {1\over 2} - 2^* \right) \sph }
\nonumber\\
&& \qquad \times
\spbk (\eqb) Z\p_\b(\lcon;k),
\label{i4.6}
\qqq
here
\qq
Z\p_\b(\lcon;K) & = &
\left( { \spk \over \spb} \right)^{N-1} \pjN Z\p_\b(\lenpqj;K),
\label{i4.7}
\\
Z\p_\b(\lenpqj;K) & = & \Leg{|\xq_j|} {\sign{\xq_j}\over \sqrt{K} }
e^{ - i \pk \kappa \spqj }
e^{ - i\pi {3\over 4} {K-2\over K} \spqj }
\label{i4.8}
\\
&&\qquad\times
\spq^{ - 4^* \Phi(U^{(p_j,q_j)} ) }
\spq^{\left(2^* - {1\over 2} \right) \spqj}
\nonumber\\
&&\qquad\times
\spajk {\spbaj\over \spk}
\spq^{4^* q_j^* (p_j \a_j^2 + s_j) }
{i\over 2}
\left( \spq^{-2^* q_j^* \a_j} - \spq^{2^* q_j^* \a_j} \right).
\nonumber\\
\qqq
In these formulas $Z\p_\b(\lenpqj;k)$ is an invariant of the
$\b$-colored link $\cL_0$ inside the lens space $\lenpqj$ constructed
by the $(\xp_j,\xq_j)$ surgery on the unknot $\cL_j$. $Z\p_\b
(\lcon;k)$ is the invariant of the $\b$-colored knot $\cL_0$ inside
the connected sum of lens spaces
\qq
\lcon = L_{-\xp_1,\xq_1} \# \cdots \# L_{-\xp_N,\xq_N}.
\label{i4.9}
\qqq
The calculation of invariants $Z\p_\b(\lenpqj;k)$
runs similar to that of $Z\p(L_{\xp,\xq};k)$ of eq.~(\ref{4.25}). The
only difference is that by taking a sum over $\a_j$ we go from
$\tU^{(\xp_j,\xq_j)}_{\a_j 1}$ to
$\tU^{(-\xq_j,\xp_j)}_{\b 1}$ rather than to
$\tU^{(-\xq_j,\xp_j)}_{1 1}$. As a result,
\qq
Z\p_\b(\lenpqj;k) & = &
\Leg{|\xp_j|} \sign{\xp_j}
\spq^{4^* \xp^*_j \xq_j - 3 s^\vee (\xq_j,\xp_j) }
\spq^{-4^* \xp_j^* \xq_j \b^2}
{\eqpb\over \eqh}
\label{i4.10}
\qqq
and
\qq
Z\p(X;k) & = &
{i\over 2\sqrt{K} }
e^{i {\pi\over 4} \kappa \sph }
e^{i\pi {3\over 4} {K-2\over K} \sph }
\spq^{\left( {1\over 2} - 2^* \right) \sph }
\label{i4.11}
\\
&&\qquad\times
\Leg{|p|} \sign{P}
\spq^{4^* P^* H - 3\sjN s^\vee(\xq_j,\xp_j) }
\nonumber\\
&&\qquad\times
{1\over \eqh}
\spbk
\qPHb
{\pjN (\eqpb) \over \eqh}
\nonumber
\qqq

The preexponential factor of eq.~(\ref{i4.11}) can be put in the form
\qq
J(\b)(\eqb),
\label{i4.12}
\qqq
where
\qq
J(\b) = \frac {\pjN (\eqpb)} { (\eqh) (\eqb)^{N-1} }.
\label{i4.13}
\qqq
It is easy to check that the function $J(\b)$ belongs to $\zq$ and
satisfies the properties of the Jones polynomial described in
Proposition~\ref{p1.2}. Therefore the full machinery of
Section~\ref{s3} could be applied to the sum of eq.~(\ref{i4.11}) in
order to convert it to the integral and ultimately prove the
Proposition~\ref{p1.1} for Seifert manifolds. However there is an
easier way. The numbers $\xp_j^*$, $\ojN$ determine a set of positive
integer numbers $\Lambda(\pps)\in \IN$ and multiplicity factors
$C_n(\pps) \in \ZZ$, $n\in \Lambda(\pps)$ such that
\qq
{\pjN (\eqpb)\over (\eqb)^{N-1} } =
\sum_{n\in \Lambda(\pps)}
C_n(\pps) (\spq^{-2^* n\b} - \spq^{2^* n\b} ).
\label{i4.14}
\qqq
Now we can apply eq.~(\ref{3.6}) to the calculation of the sum
\qq
\lefteqn{
{1\over \eqh} \spbk \qPHb {\pjN(\eqpb)\over (\eqb)^{N-2} }
}
\label{i4.15}
\\
&=&
{2\over\eqh} \sum_{n\in\Lambda} C_n \spbk \qPHb
(\spq^{-2^* (n+1) \b} - \spq^{-2^* (n-1) \b})
\nonumber\\
&=&
2\sqrt{K} \Leg{P^* H} \eikapw \sum_{n \in \Lambda} C_n
\frac{\spq^{4^* PH^* (n+1)^2} - \spq^{4^* PH^* (n-1)^2} } {\eqh}
\nonumber\\
&=&
2\sqrt{K} \Leg{P^* H} \eikapw \sum_{n \in \Lambda} C_n
\spq^{4^* PH^* (n^2 + 1)}
\frac{\spq^{2^* PH^* n} - \spq^{-2^* PH^* n} } {\eqh}.
\nonumber
\qqq
Combining eqs.~(\ref{i4.11}) and~(\ref{i4.15}) we find the formula
\qq
Z\p(X;k) & = &
\Leg{|H|} \sign{H} \lqexp
\label{i4.16}
\\
&&\qquad\times
\sum_{n\in \Lambda} C_n q^{4^* PH^* (n^2+1)}
\frac{ \spq^{-2^* PH^* n} - \spq^{2^* PH^* n}} {\eqh},
\nonumber
\qqq
which demonstrates that $Z\p(X;k)\in \zq$.

Now we come back to eqs.~(\ref{i4.14}),~(\ref{i4.15}) and use the
fact that for $n\in \ZZ$,
\qq
\spbk \spq^{-4^* P^* H \b^2 - 2^* n\b} & = &
\sqrt{K} \Leg{P^* H} \eikapw \spq^{4^* PH^* n^2},
\label{i4.17}\\
\intinfb \spq^{-{\b^2\over 4H^* P} - 2^*\b n} & = &
e^{ - i{\pi\over 4} \sign{PH^*} }
\sqrt{2K|PH^*|} \spq^{4^* PH^* n^2},
\label{i4.18}
\qqq
so that
\qq
\spbk \spq^{-4^* P^* H \b^2 - 2^* n\b} =
\Leg{P^* H} \eikapw \frac{ e^{i {\pi\over 4} \sign{PH^*} } }
{\sqrt{2|PH^*|} }
\intinfb \spq^{-{\b^2\over 4H^* P} - 2^*\b n}.
\label{i4.19}
\qqq
This equation allows us to convert the sum over $\b$ in
eq.~(\ref{i4.15}) into an integral. Then by using eq.~(\ref{i4.14})
backwards we arrive at eq.~(\ref{i4.11}) with the integral instead of
a sum:
\qq
Z\p(X;k) & = &
\Leg{|H|} \sign{H} \lqexp I(\spq),
\label{i4.21}
\\
I(\spq) & = &
{1\over \eqhm} {e^{i{\pi\over 4}\sign{PH^*} }\over 2\sqrt{2K|PH^*|}}
\intinfb \spq^{- {\b^2\over 4PH^*} }
{\pjN (\eqpb)\over (\eqb)^{N-2} }.
\label{i4.22}
\qqq
The integral over $\b$ is well defined in view of eq.~(\ref{i4.16})
(actually, one might add a regularizing factor
$\lim_{\epsilon\rightarrow 0} e^{-\epsilon\b^2}$. It can also be
calculated by expanding the preexponential factor in powers of
$x=q-1$ and integrating their coefficients, which are polynomials in
$\b$, with the gaussian exponential
$\spq^{- {\b^2\over 4PH^*}}$. This procedure leads to the following
relation:
\qq
[I(\spq)]^\diamond & = &\left[
{1\over\eqhm} {e^{i{\pi\over 4}\sph} \over 2\sqrt{2K} }
\sqrt{\left| {H\over P}\right| }
\spint{\b} d\b\,
\spq^{-{1\over 4} {H\over P} \b^2 }
\frac {\pjN \left( \spq^{-{\b\over 2p_j}} - \spq^{\b\over 2p_j}
\right)}
{\left( \spq^{-{\b\over 2}} - \spq^{\b\over 2} \right) ^{N-2} }
\right]^\vee
\label{i4.23}
\\
& = &
\left[ {1\over \sqrt{2\over K} \spk}
{\sph\over K}
e^{i {3\over 4} \pi \sph }
\sqrt{\left| {H\over P}\right| }
\lint
\right]^\vee
{}.
\nonumber
\qqq
Since
\qq
\left[ \lqexp \right]^\diamond
=
\left[
\spq^{ {1\over 4} {H\over P} - {3\over 4} \sph -
3\sjN s(\xq_j,\xp_j)}
\right]^\vee
\label{i4.24}
\qqq
and according to ~\cite{Ro*1},~\cite{Ro*2},~\cite{Ro1} (see, e.g.
eq.(4.9) of~\cite{Ro1})
\qq
\ztr(X;k) & = & {e^{i\pi{3\over 4} \sph} \over K}
{\sign{P}\over \sqrt{|P|} }
e^{{i\pi\over 2K} \left[ {H\over P} - 3\sph - 12 \sjN
s(\xq_j,\xp_j)\right] }
\nonumber\\
&&
\qquad\times
\lint,
\label{i4.25}
\qqq
we conclude that eq.~(\ref{1.21}) holds for Seifert manifolds which
are rational homology spheres.

\nsection{Discussion}
\label{s5}

As we have already mentioned in the Introduction, all results of this
paper follow rigorously from the ``physical'' input:
Propositions~\ref{p1.2} and~\ref{p1.3}. The Proposition~\ref{p1.3}
might be harder to prove with mathematical rigor. Even its
formulation uses the asymptotic structure~(\ref{1.15}) of Witten's
invariant at large values of $K$. To our knowledge, this structure
has not been rigorously derived yet from the surgery
formula~(\ref{1.2}). The Proposition~\ref{p1.2} follows from the
properties of Reshetikhin's formula~~\cite{Ro2} of the Jones
polynomial of a link and seems to have better chances for
legitimization.

If we forget for a moment about the trivial connection contribution
to Witten's invariant, then we can use eq.~(\ref{1.25}) as a
definition of the manifold invariants $\snm$. The
Proposition~\ref{p1.2} gives us enough information in order to derive
eq.~(\ref{1.21}) in the form
\qq
\lefteqn{
\left[ \ordH \Leg{\ordH} \zpmk \right]^\diamond
}
\hspace*{2in}
\label{5.1}
\\
&=&
\left[ \ordH^{3\over 2} \exp \left( (\snm - \snst) \pipk^n
\right) \right]^\vee.
\nonumber
\qqq
This also proves Ohtsuki's Theorem~\ref{t1.2}, since we can define
the invariants $\lnm$ by eq.~(\ref{1.18}). According to the comments
of~\cite{Oh1},~\cite{Oh2}, we may also conclude that $\snm$, as
defined by the surgery formula ~(\ref{1.25}), are indeed invariants
of $M$.

In our previous paper~\cite{Ro4} we studied the properties of $\snm$
as they follow from eq.~(\ref{1.25}), Proposition~\ref{p1.2} and some
other properties of Reshetikhin's formula~\cite{Ro2}. Now we use
eq.~(\ref{1.18}) in order to extend them to $\lnm$:
\begin{proposition}
\label{p5.1}
The invariants $\lnm$ are finite type invariants of \Rhs as defined
in~\cite{Oh3} and~\cite{Ga2}. An invariant $\lnm$ is of Ohtsuki order
$3n$, \opr (~\cite{Ro4}) order $2n$ and at most of Garoufalidis
order $n$. Also
\qq
&\displaystyle
2^{4n}\,n!\,(2n)!\,(9n)!\, \ordH^n \lnm \in \ZZ,
\label{5.2}
\\
&\displaystyle
\ordH^n \lnm\in \ZZ\left[ {1\over 2},{1\over 3},\ldots,{1\over 2n}
\right].
\label{5.3}
\qqq
\end{proposition}

To see that $\lnm$ is of exactly Ohtsuki order $3n$ one has to
find an $n$-component link $\cL$ such that the alternating sum
\qq
\sum_{\cL\p\subset \cL} (-1)^{\#\cL\p}
\lambda_n(\chi_{\cL\p} (S^3) )
\label{5.4}
\qqq
is non-zero (here $\chi_{\cL}(S^3)$ denotes a manifold (\Rhs)
constructed by $(1,1)$ surgeries on the components of a link $\cL$
in $S^3$). Recall that according to eq.~(\ref{1.18})
\qq
\lnm =
\sum_{m_1,\ldots,m_n \geq 0
\atop m_1 + 2m_2 + \cdots + nm_n = n}
C^n_{m_1,\ldots,m_n} S_1^{m_1}(M) \cdots S_n^{m_n}(M)
\label{5.5}
\qqq
for some numbers  $C_{m_1,\ldots,m_n}$. Consider the $n$-loop diagram
consisting of $(n-1)$ small loops sitting on one big loop, and
the corresponding~\cite{Oh3} link $\cL$. This diagram has no
subdiagrams with only trivalent vertices. Then according
to~\cite{Ro4}, for any $S_{n\p}(M)$, $n\p<n$ the alternating sums are
equal to zero
\qq
\sum_{\cL^{\prime\prime} \subset \cL\p} (-1)^{\#\cL^{\prime\prime}}
S_{n\p} (\chi_{\cL^{\prime\prime}}(M) ) = 0,
\label{5.6}
\qqq
if $\cL\p\subset L$ and $\#\cL\p \geq n\p$. Therefore an invariant
$S_{n\p}(M)$ is of Ohtsuki order less than $n\p$ with respect to
$\cL$ and its subdiagrams (of course, there exist other links for
which the \lhs of eq.~(\ref{5.6}) is non-zero if $\#\cL\p = n\p$). As
a result, only the term $C_{0,\ldots,0,1} S_n(M)$ in eq.~(\ref{5.5})
matters for the calculation of the alternating sum~(\ref{5.4}).
Since, according to~~\cite{Ro4}, the sum
\qq
\sum_{\cL\p\subset \cL} (-1)^{\#\cL\p}
S_n(\chi_{\cL\p} (S^3) )
\label{5.7}
\qqq
is non-zero, we conclude that the sum~(\ref{5.4}) is also non-zero.
Hence $\lnm$ is of exactly Ohtsuki order $3n$.

Finally, let us comment on the relation~(\ref{1.21}) between the
invariant $\zpmk$ at prime values of $K$ and the trivial connection
contribution to Witten's invariant. This relation is not an obvious
result in the sense that the contributions of non-trivial connections
to $\zpmk$ do not seem to cancel at prime $K$ (one might expect that
there is some cancelation leaving only the contribution of the
trivial connection). This can be seen at the example of a lens space
$L_{\xp,\xq}$ for which
\qq
Z\p(L_{\xp,\xq})
 & = & \sign{\xp} \Leg{|\xp|}
\,e^{{6\pi i\over K} s^\vee(\xq,\xp)}
\,
(-1)^{\xp^* + 1} \,
\frac{\sin\left( {\pi\over K} \xp^*\right)}
{\sin\left({\pi\over K}\right)},
\label{5.8}
\\
\frac {Z^{({\rm tr})} (L_{\xp,\xq};k)}
{Z(S^3;k)}
& = &
{\sign{\xp}\over \sqrt{|p|}}
\,e^{{6\pi i\over K} s(\xq,\xp)}
\frac{\sin\left( {\pi\over K} {1\over p}\right)}
{\sin\left({\pi\over K}\right)}.
\label{5.9}
\qqq
Although these two expressions have many common features, which
ultimately lead to the relation~(\ref{1.21}), still their numerical
values are quite different. The $p$ and $K$ dependence of
$Z\p(L_{\xp,\xq})$ is somewhat typical of contributions of
$U(1)$-reducible connections. Besides, the Dedekind sum $s(\xq,\xp)$
is generally a fraction, so $s^\vee(\xq,\xp)\neq s(\xq,\xp)$.

We established the relation~(\ref{1.21}) by comparing directly the
surgery formulas. It would be much better to have a conceptual
explanation for this phenomenon. One might speculate that it would
come from number theory and perhaps $p$-adic quantum field theories.

\section*{Acknowledgements}

I am very thankful to J.~Roberts for many stimulating discussions and
for attracting my attention to the results of~\cite{Mu1}
and~\cite{Mu2}. I also appreciate the numerous conversations with
D.~Freed, P.~Freund, C.~Gordon, L.~Kauffman, N.~Reshetikhin,
A.~Vaintrob, O.~Viro, K.~Walker and their comments. I am especially
thankful to O.~Alvarez, L.~Mezincescu and R.~Nepomechie for their
constant encouragement and friendly support.

This work was supported by the National Science Foundation
under Grant No. PHY-92 09978.

\end{document}